\documentclass{article}
\usepackage{natbib}
\usepackage{graphicx}

\usepackage{amsmath,amsfonts,amsthm}
\newtheorem{lemma}{Lemma}[section]
\newtheorem{proposition}{Proposition}[section]
\newtheorem{definition}[proposition]{Definition}
\newtheorem{theorem}[proposition]{Theorem}
\newtheorem{remark}[proposition]{Remark}
\newcommand{\DT}{\mathrm{\Delta}t}
\newcommand{\pder}[2]{\frac{\partial #1}{\partial #2}}
\newcommand{\diag}{\mathsf{diag}}
\newcommand{\tridiag}{\mathsf{tridiag}}

\newcommand{\R}{\mathbb{R}}
\newcommand{\C}{\mathbb{C}}
\newcommand{\imm}{\mathrm{i}}

\newcommand{\tr}{\mathsf{T}}

\renewcommand{\vec}[1]{\mathbf{#1}}

\title{Multigrid and preconditioning strategies for implicit PDE
  solvers for degenerate parabolic equations}

\author{%
{\sc M. Donatelli\thanks{Email: marco.donatelli@uninsubria.it},
M. Semplice\thanks{Corresponding author. Email:
  matteo.semplice@uninsubria.it}
and
S. Serra-Capizzano\thanks{Email: stefano.serrac@uninsubria.it}
}\\[2pt]
Dipartimento di Fisica e Matematica\\ Universit\`a
    dell'Insubria \\ Via Valleggio 11, 22100 Como, Italy.
}

\begin{document}

\maketitle

\begin{abstract}
  {%
    The novel contribution of this paper relies in the proposal of a
    fully implicit numerical method designed for nonlinear degenerate
    parabolic equations, in its convergence/stability analysis, and in
    the study of the related computational cost. In fact, due to the
    nonlinear nature of the underlying mathematical model, the use of
    a fixed point scheme is required and every step implies the
    solution of large, locally structured, linear systems. A special
    effort is devoted to the spectral analysis of the relevant
    matrices and to the design of appropriate iterative or
    multi-iterative solvers, with special attention to preconditioned
    Krylov methods and to multigrid procedures: in particular we
    investigate the mutual benefit of combining in various ways
    suitable preconditioners with V-cycle algorithms.  Numerical
    experiments in one and two spatial dimensions for the validation
    of our multi-facet analysis complement this contribution.}
{\bf AMS SC: 65N12, 65F10  (65N22, 15A18, 47B35)}
\end{abstract}

\section{Introduction}

We  consider a single equation of the form
\begin{equation}\label{eq:parab}
\pder{u}{t} = \nabla\cdot\left(D(u)\nabla u\right),
\end{equation}
where $D(u)$ is a non-negative function. The the equation is parabolic
and it called
degenerate whenever $D(u)$ vanishes for some values of $u$.
 For the convergence analysis of
our numerical methods, we will require that $D(u)$ is at least
differentiable and that $D'(u)$ is Lipschitz continuous, while the
existence of solutions is guaranteed under the milder assumption of
continuity \citep[see][]{VazquezBOOK}.

The classical porous medium equation (where $D(u)$ is restricted to be
a power law) or its generalized form have important applications in
many fields of science. Their name arise from the model of Darcy's
flow of a gas through a porous medium ($D(u)=\kappa u^\gamma$ where
$\gamma$ is the specific heat ratio), but classical applications range
from ground-water flow (Boussinesq equation) to spatial population
dynamics (where crowding effects require nonlinear degenerate
diffusion terms). Moreover models based on equation \eqref{eq:parab}
have been proposed as useful approximations of more complex models
like thin films motion (when disregarding surface tension), water-oil
mixtures in porous medium, boundary layer in fluid flow past an
obstacle, magma in volcanoes, etc. More recently they have been
studied as limits of kinetic particle models, with applications for
the diffusion in semiconductors. Finally we mention that some
contrast-enhancement filters for image processing like the one by
Perona-Malik are based on \eqref{eq:parab}. For more details on the
aplications, see e.g.  \citet[Chap 2 and 21]{VazquezBOOK} and the
references therein.

The present investigation is part of the search for suitable numerical
techniques to integrate for long times nonlinear, possibly degenerate,
parabolic equations appearing in models for monument degradation (see
\citet{ADN:sulfation}) when chemical/micro-biological pollutants are
taken into consideration. We wish to point out that the techniques
developed here have applications that go beyond the aforementioned
models. For example, again in the area of conservation of the cultural
heritage, they could be adapted to numerically investigate the more
complete sulfation model described in \citet{AFNT07:model} and the
consolidation model presented in \citet{CGNNS:teos}. Some applications
in the field of monument conservation have been presented in
\citet{SDS:monum}, where the mathematical tools developed in the
present paper are employed for forecasting marble deterioration.  Of
course, in such a context, given the wide variety of artefacts, an
important challenge is the combination of the approximation scheme
with the related linear algebra solvers, in presence of complicate
geometries and griddings.

In the literature, degenerate parabolic equations have been
discretized mainly using explicit or semi-implicit methods, thus
avoiding to solve the nonlinear equation arising from the elliptic
operator. A remarkable class of methods arise directly from the
so-called non-linear Chernoff formula \citet{BP72} for time
advancement, coupling it with a spatial discretization: for finite
differences the latter study was started in \citet{BBR79} and for
finite elements in \citet{MNV87}. More recently, another class related
to the relaxation approximation emerged: such numerical procedures
exploit high order non-oscillatory methods typical of the
discretization of conservation laws and their convergence can be
proved making use of semigroup arguments similar to those relevant for
proving the Chernoff formula \citep{CNPS07:degdiff}.

In this paper we start from the Crandall-Liggett formula
\begin{equation}\label{CLformula}
U(t^n,x)-\DT L_D(U(t^n,x)) = U(t^{n-1},x),
\end{equation}
where time has been discretized with steps $\DT=t^n-t^{n-1}$, and
$-L_D(\cdot)$ denotes the elliptic operator $u\mapsto
-\nabla\cdot(D(u)\nabla u)$. The computation of the numerical
solution $U(t^n,x)$ requires to solve a nonlinear equation whose
form is determined by the elliptic operator and the nonlinear
function $D(u)$, but the convergence is guaranteed without
restrictions on the time step $\DT$ \citep{CL71}.  Furthermore, due
to the nonlinear nature of the underlying mathematical model, the
use of a fixed point scheme is required and the choice of the
Newton-like methods implies the solution at every step of large,
locally structured \citep[in the sense of][]{Tilli98} linear
systems. A special effort is devoted to the spectral analysis of the
relevant matrices and to the design of appropriate iterative or
multi-iterative solvers \citep[see][]{Serra93:multi}, with special
attention to preconditioned Krylov methods and to multigrid
procedures (see
\citet{Greenbaum:book,Saad:book,Hackbusch:book,Trottenberg:book} and
references therein for a general treatment of iterative solvers).
Although most of the analysis is developed in the one-dimensional
case (from Section \ref{sec:FD} to Section \ref{sec:NUM}), we also
indicate in Section \ref{sec:2D} how to generalize our approach to
two spatial dimensions. we also perform numerical experiments for the
validation of our analysis in both settings (see Section
\ref{sec:NUM} and Section \ref{sec:2D}, respectively).

The paper is organized as follows. In Section \ref{sec:FD} we couple
the time discretization \eqref{CLformula} with a spatial
discretization based on finite differences and set up a Newton method
for the resulting system of nonlinear equations.  We report the
explicit form of the Jacobian appearing in the Newton iterations and
we prove the convergence of the Newton methods under a mild
restriction on $\DT$.  In Section \ref{sec:LIN} we consider various
iterative methods for the solution of the inner linear systems
involved in the Newton method. A brief spectral analysis of the
related matrix structures is provided in order to give an appropriate
motivation for the good behaviour of the proposed iterative and
multi-iterative solvers. In Section \ref{sec:NUM} we perform some
numerical tests. In Section \ref{sec:2D} we describe a generalization
of the previous methods to a two-dimensional case and perform
numerical tests in this setting too.  Finally, a conclusion section
with a short plan for future investigations completes the paper.

\section{Numerical methods in one dimension}
\label{sec:FD} 
In order to discretize equations like \eqref{eq:parab}, we will employ
a time semi-discretization given by the Crandall-Liggett formula and a
space discretization based on finite differences, explained in the
following subsection. The latter numerical choice leads to a system of
coupled nonlinear equations that need to be solved at each discrete
timestep in order to compute the solution of the PDE: this is achieved
using the Newton method, as detailed in Subsection \ref{ssec:newton},
where we also prove and comment convergence results.

\subsection{Finite difference discretization}\label{ssec:FD}
We take into consideration a standard discretization in space using
finite differences. Denoting $x_\xi = a+\xi h$, we consider $N+2$
points with equal spacing $h=(b-a)/(N+1)$ in the interval $[a,b]$ and
we denote by $u^n_k$ the approximate solution at time $t^n$ and
location $x_k$, where $k=0,\ldots,N+1$. Let $\vec{u}^n$ be the vector
containing the collection of the unknown values $u^n_k$. When no
potential confusion can arise, we will sometimes drop in both
notations the superscript indicating the time level. Of course, when
considering Dirichlet boundary conditions, the values $u_0$ and
$u_{N+1}$ are known and can be eliminated by the equations, leaving a
vector of unknowns $\vec{u}^n$ that contains only $u_k^n$ for
$k=1,\ldots,N$.  Boundary conditions of Neumann or Robin type can be
treated in similar ways.

We choose a standard 3-points second order approximation of the
differential operator $(D(u)u_x)_x$. Denoting with the subscript
$\xi$ the evaluation at the point $x_\xi$, we have that:
\begin{multline}\label{eq:FDlaplacian}
\left.\pder{}{x}\left(D(u)\tfrac{\partial u}{\partial x}
\right)\right|_j
= \frac{D(u)_{j+1/2} \left.\tfrac{\partial u}{\partial x}\right|_{j+1/2}
-D(u)_{j-1/2} \left.\tfrac{\partial u}{\partial x}\right|_{j-1/2}}{h}+o(1)\\
=\frac{D(u)_{j+1/2} (u_{j+1}-u_j)
-D(u)_{j-1/2} (u_j-u_{j-1})}{h^2}+o(1)\\
=\frac{(D(u_{j+1})+D(u_j)) (u_{j+1}-u_j) -(D(u_j)+D(u_{j-1}))
(u_j-u_{j-1})}{2h^2}+o(1)
\end{multline}
where the $o(1)$ error term is of order $h^2$ under the assumption
that the composition
\[\phi(\cdot)=D(u(\cdot))\]
is at least continuously differentiable, with Lipschitz continuous
first derivative. Putting together all the contributions for different
grid points, we end up with $L_{D(\vec{u})}\vec{u}$, where the tridiagonal matrix
\begin{equation}\label{LD}
L_{D(\vec{u})} =
\begin{bmatrix}
-D_{1/2}-D_{3/2} & D_{3/2} & \\
D_{3/2} & -D_{3/2}-D_{5/2} & D_{5/2} & \\
 &   D_{5/2} &  \ddots & \ddots & \\
&&\ddots & \ddots & D_{N-1/2} \\
 & & & D_{N-1/2} & -D_{N-1/2} -D_{N+1/2} \\
\end{bmatrix}
\end{equation}
conatins the values
\begin{equation*}
D_{j+1/2}=\frac{D(u_{j+1})+D(u_{j})}2
\,,\qquad j=0,\ldots,N,
\end{equation*}
and thus depends nonlinearly on the $u_j$'s.  It should be noticed
that the latter is a second order approximation of $\phi(x_{j+1/2})$
since $u^n_{k}$ differs from $u(t^n,x_k)$ by $O(h^2)$ thanks to the
second order scheme and since, by standard Taylor expansions, we have
\begin{eqnarray*}
D_{j+1/2}=\frac{D(u_{j+1})+D(u_{j})}2 & = &
\frac{\phi(x_{j+1})+\phi(x_{j})}2+O(h^2)\\
& = & \phi(x_{j+1/2})+ \frac{h^2}{8}\phi_{xx}(\eta(h,j))+O(h^2) \\
& = & \phi(x_{j+1/2})+ O(h^2),\ \quad \eta(h,j)\in (x_{j-1},x_{j}),
\end{eqnarray*}
under the mild assumption that $\phi_{xx}(\cdot)$ is a bounded
function. Of course, the same conclusion holds if $\phi_{x}(\cdot)$
is Lipschitz continuous.

In the following, we denote by $\tridiag_k[\beta_k, \alpha_k,
\gamma_k]$ a square tridiagonal matrix of order $N$ with entries
$\beta_k$ on the lower diagonal, $k=2,\cdots,N$, $\alpha_k$ on the
main diagonal, $k=1,\cdots,N$, and $\gamma_k$ on the upper diagonal,
$k=1,\cdots,N-1$. With this notation, $L_{D(\vec{u})} =\tridiag_k
[D_{k-1/2}, -D_{k-1/2}-D_{k+1/2}, D_{k+1/2}]$. We also denote with
$\diag_k[\alpha_k]$ the square diagonal matrix with $\alpha_k$ on the
$k^{\text{th}}$ row.

As already observed $L_{D(\vec{u})}$ is a symmetric real tridiagonal
matrix. Since $D(\cdot)$ is a nonnegative function, the matrix
$-L_{D(\vec{u})}$ is always positive semidefinite, beacuse it is
weakly diagonally dominant by row, or equivalently thanks to the first
Gerschgorin Theorem \citep[see e.g.][]{Golub:book}. Furthermore, we
have positive definiteness (i.e. invertibility), at least for every
$N$ large enough, if in addition $\phi(\cdot)$ has only isolated zeros
in $(a,b)$. In that case, for $N$ large enough, the matrix is
irreducible or block diagonal with irreducible blocks.  In particular,
when $\phi(\cdot)$ is strictly positive in $(a,b)$ then
$-L_{D(\vec{u})}$ is positive definite and irreducible for any $N$.

When introducing numerical methods for the approximation of the
differential equation we will encounter nonlinear systems involving
the matrices $-L_{D(\vec{u})}$. At that point more sophisticated
(spectral) relations and features will be discussed, when choosing the
appropriate iterative solvers for the global linearised system (see
Subsection \ref{ssec:SPECTRAL}). For the moment we just observe that,
thanks to the previous preliminary spectral analysis, all the
classical iterative solvers like Jacobi and Gauss-Seidel (and their
damped version with damping parameter belonging to $(0,2)$) are all
convergent for the solution of a linear system with such a coefficient
matrix. The problem is that the spectral radii are very close to $1$,
with a gap ranging between $O(N^{-2})$, reached by all these classical
iterations with the only exception of the optimally damped
Gauss-Seidel, and $O(N^{-1})$, reached for Gauss-Seidel with optimal
damping parameter \citep[see][]{Varga:matrixbook}. When considering
the whole system things become slightly better since the gap between
the spectral radius and $1$ reduces for all the considered procedures
to $O(N^{-1})$.  However, as a partial conclusion, we can safely claim
the considered iterations would be unacceptably slow and the search
for specialised iterative solvers becomes mandatory. This latter is the
main subject of Section~\ref{sec:LIN}.

\subsection{The nonlinear system and the Newton iteration}
\label{ssec:newton}

Following the Crandall-Liggett formula \eqref{CLformula}, in order
to compute $\vec{u}^n$ from $\vec{u}^{n-1}$, we need to solve the
nonlinear vector equation
\begin{equation*}
\vec{u}^n =\vec{u}^{n-1} + \frac{\DT}{h^2} L_{D(\vec{u}^n)}\vec{u}^n
\end{equation*}
and thus we set up Newton iterations for the vector function
\begin{equation}\label{eq:FnonlinD}
F(\vec{u}) = \vec{u} -\frac{\DT}{h^2} L_{D(\vec{u})}\vec{u} -
\vec{u}^{n-1}.
\end{equation}
In the following, we denote $\vec{u}^{n,s}$ the $s^{\text{th}}$
Newton iterate for the computation of $\vec{u}^{n}$.
The generic partial derivative of $F(\vec{u})$ is
\begin{equation} \label{eq:JacD}
  \pder{F_k}{u_j} 
  =\delta_{jk}-\frac{\DT}{h^2}\left.L_{D(\vec{u})}\right|_{j,k}
  -\frac{\DT}{2h^2}\left[
    \begin{aligned}
    &\delta_{k-1,j}D^{\prime}_{k-1}(u_{k-1}-u_k)+\\
    &+\delta_{k,j}D^{\prime}_{k}(u_{k-1}-2u_k+u_{k+1})+\\
    &+\delta_{k+1,j}D^{\prime}_{k+1}(u_{k+1}-u_k).
  \end{aligned}
\right]
\,,
\end{equation}
so that the Jacobian is
\begin{align}
F^\prime(\vec{u}) &= X_N(\vec{u}) + Y_N(\vec{u}),\label{Jf}
\\
X_N(\vec{u}) &= I_N-\frac{\DT}{h^2}L_{D(\vec{u})}, \label{Xn}
\\
 Y_N(\vec{u}) &= -\frac{\DT}{2h^2}T_N(\vec{u}) \diag_k[D^\prime_k],\label{Yn}
\\
T_N(\vec{u}) &=
\tridiag_k[u_{k-1}-u_k,u_{k-1}-2u_k+u_{k+1},u_{k+1}-u_k].\label{Tn}
\end{align}
The matrix $X_N(\vec{u})$ is symmetric positive definite and
$\lambda_{\text{min}}(X_N(\vec{u}))\geq1$, where
$\lambda_{\text{min}}(A)$ denotes the minimum eigenvalue of the
matrix $A$. We note that the inequality is strict under assumption
of isolated zeros. If $\diag_k[D^\prime_k]$ is positive
semidefinite, i.e. if $D(\cdot)$ is a smooth nondecreasing
function, then, setting $E^2=\diag_k[D^\prime_k]$, $E$ is a positive
semidefinite diagonal matrix and $Y_N(\vec{u})$ is similar to
$-\frac{\DT}{2h^2}ET_N(\vec{u})E$. Moreover, defining
\begin{align}
 \tilde Y_N(\vec{u}) &= 
 -\frac{\DT}{2h^2}\tilde T_N(\vec{u}) \diag_k[D^\prime_k],\label{tildeYn}
\\
\tilde T_N(\vec{u}) &=
\tridiag_k[u_{k-1}-u_k,0,u_{k+1}-u_k]=T_N(\vec{u})
-\diag_k[u_{k-1}-2u_k+u_{k+1}],\label{tildeTn}
\end{align}
we have that $\tilde Y_N(\vec{u})$ is similar to $-\frac{\DT}{2h^2}
E \tilde T_N(\vec{u})E$, with the latter being anti-symmetric, which
implies a pure imaginary spectrum.

In the following we will denote by $\|\cdot\|$ the Euclidean norm
for vectors and the induced spectral norm for matrices.

\begin{remark} \label{remark:Yn}
  If $\vec{u}$ is a sampling of a solution $u$ of \eqref{eq:parab}
  at least continuous and $\omega_u(\cdot)$ denotes its modulus
  of continuity, then
  \[
 \| Y_N(\vec{u}) \| \le e_u(\DT,h)
  \]
  with
  \[
  e_u(\DT,h)=4\frac{\DT}{h^2}\|D'(u)\|_\infty\omega_u(h).
  \]
  In order to deduce the latter, it is enough to recall that for normal
  matrices the spectral norm is bounded by any induced
  norm and in particular by the the matrix norm induced by the
  infinity vector norm. In particular,
  if $u$ is H\"older continuous with exponent $\alpha\in (0,1]$ and
  constant $M>0$, then the estimate above can be written as
  \[
  e_u(\DT,h)=4 M\frac{\DT}{h^{2-\alpha}}\|D'(u)\|_\infty.
  \]
  Of course if $u$ is continuously differentiable,
  then we find $\alpha=1$ and $M=\|u'\|_\infty$. Furthermore, when $u$
  is two times continuously differentiable, we notice that
  $u_{k-1}-2u_k+u_{k+1}=h^2 u''(\xi_k)$ which leads to a more refined
  expression i.e.
  \[
  e_u(\DT,h)=\left(2 \frac{\DT}{h}\|u'\|_\infty+\DT \|u''\|_\infty\right)\|D'(u)\|_\infty.
  \]
  Finally, if we are interested in evaluating $\|Y_N(\vec{\tilde
    u})\|$ where $\tilde u$ is an approximation to the true solution
  $u$ (this happens naturally in the numerical process discussed in
  the present section), then
  \[
  \|Y_N(\vec{\tilde u})\|\le e_{\tilde u}(\DT,h)\le
  e_u(\DT,h)+4 \frac{\DT}{h^2}\|D'(u)\|_\infty \|u-\tilde u\|_\infty
  = 4 \frac{\DT}{h^2}\|D'(u)\|_\infty\left(\omega_u(h)+\|u-\tilde
  u\|_\infty\right).
  \]
  Hence, since we are using second order formulae,
  the error $\|u-\tilde u\|_\infty=O(h^2)$ and therefore $\|Y_N(\vec{\tilde u})\|$
  is dominated by $\omega_u(h)$, which is of order $h$
  if the solution is Lipschitz continuous, that is
  \[
  \|Y_N(\vec{\tilde u})\|\le 4 M\frac{\DT}{h}\|D'(u)\|_\infty + O(\DT).
  \]
  In conclusion, we can safely claim that the global spectrum of the
  Jacobian $F^\prime(\vec{\tilde u})$ is decided, up to small perturbations,
  by the matrix $X_N(\vec{\tilde u})$. For making more explicit the latter
  statement, if we assume that $\DT=C h$, where $C>0$ is independent of
  $h$, then $\lambda_{\text{min}}(X_N(\vec{u}))\geq1$,
  $\|X_N(\vec{\tilde u})\|= O(h^{-1})$, while
  $\|Y_N(\vec{\tilde u})\|=O(1)$.
\end{remark}

In order to prove the convergence of the Newton method, we first
consider some auxiliary results.

\begin{lemma} \label{lem:sigmaA}
For a generic matrix $A$, the minimum singular value is
\begin{equation}\label{eq:smin}
    \sigma_{\rm min}(A) \geq \lambda_{\rm min}\left(\frac{A+A^\tr}{2}\right).
\end{equation}
\end{lemma}
\begin{proof}
Consider the symmetric matrix
\[B=\left[
    \begin{aligned}
    & 0 & A^\tr\\
    & A & 0
  \end{aligned}
\right]\]
with eigenvalues $\lambda_1(B)=\sigma_1(A) \geq \dots \geq
\lambda_N(B)=\sigma_N(A)\ \geq \lambda_{N+1}(B)=-\sigma_N(A)\geq
\dots \geq \lambda_{2N}(B)=-\sigma_1(A)$, since the Schur
decomposition of $B$ is easily written in terms of the singular
value decomposition of $A$ \citep[see][]{Bhatia:book,Golub:book}.
Let $\cal V$ be a vector space and $\vec{x} \in \R^N$ such that
$\|\vec{x}\|_2>0$. Thanks to the minimax principle
\citep{Bhatia:book}, we obtain
\[
\sigma_{\rm min}(A)  =  \lambda_N(B) = \max_{\dim({\cal V})=N} \;
\min_{\vec{y}\in {\cal V}}
\frac{\vec{y}^TB\vec{y}}{\vec{y}^T\vec{y}}
 \geq  \min_{\substack{\vec{y} = [\vec{x}^T\!\!, \; \vec{x}^T]^T\\\vec{x}\in\mathbb{R}^N}}
    \frac{\vec{y}^TB\vec{y}}{\vec{y}^T\vec{y}}
     = \min_{\vec{x}\in \R^{N}}\frac{\vec{x}^T(A^T+A)\vec{x}}{2\vec{x}^T\vec{x}}
     =  \lambda_{\rm min}\left(\frac{A+A^T}{2}\right).
  \]
\end{proof}

\begin{remark}\label{tools-algebra}
The proof technique used for bounding from below the minimal
singular value of a matrix $A$ is part of a more general
framework useful for refining, when necessary, the estimates. In
fact, in general it can be proved that for any complex-valued matrix
$A$ the minimal singular value is not less than the distance $d_r$
of any straight line $r$ separating the numerical range of $A$ from
the complex zero. Therefore a better estimate can be obtained by
computing the sup (that we call $d$) of $d_r$, over all straight lines
that induce the separation. In our case we used the fact that the
real part of $A$ (that is ${\rm Re}(A)=(A+A^\tr)/2$) is positive
definite and so our straight line becomes the set of all complex
numbers having real part equal to $\lambda_{\rm min}({\rm Re}(A))$.
The estimate could be poor since the latter straight line is not
necessarily tangent to the numerical range (a convex set by the
Toeplitz-Hausdorff theorem, see \citet{Bhatia:book}): thus $d$ could be
much larger than $d_r$. However in our setting such an estimate is
already very satisfactory, as also stressed by the numerical
experiments.
\end{remark}

\begin{proposition}\label{prop:Fprimo}
  Consider $F(\vec{u})$ as defined in \eqref{eq:FnonlinD}, where
  $\vec{u}$ is a sampling (at a given time $t$) of a solution $u$ of
  \eqref{eq:parab} with $D$ differentiable and having first derivative
  Lipschitz continuous.  If, in addition, $\vec{u}$ is differentiable
  with Lipschitz continuous first derivative, then
  \begin{equation}\label{eq:Finvnorm2}
    \left\| F^\prime(\vec{u})^{-1}\right\| \leq 1+O(\DT).
  \end{equation}
    When using the induced $l^\infty$ norm, we have
  \begin{equation}\label{eq:Finvnorminf}
    \left\|
    F^\prime(\vec{u})^{-1}\right\|_\infty \leq C_1
  \end{equation}
  for $h$ sufficiently small and under the additional assumption that $\DT\le C_\infty h$ for some $C_\infty >0$.
\end{proposition}
\begin{proof}
  For the sake of notational simplicity, we set
  $A=F^{\prime}(\vec{u})$.  First of all we write the symmetric
  part of $A$ as
  \[ \frac{A+A^\tr}{2} = X_N(\vec{u})+Z_N(\vec{u}) \,,\]
  where
  \[
  Z_N(\vec{u}) = - \frac{\DT}{4h^2} \, \tridiag_k\left[
      (D^{\prime}_{k-1}-D^{\prime}_k)(u_{k-1}-u_k),
      2D^{\prime}_{k}(u_{k-1}-2u_k+u_{k+1}),
      (D^{\prime}_{k+1}-D^{\prime}_k)(u_{k+1}-u_k)
  \right].
  \]
  By the regularity of $D$ and $u$,
  we have that every entry of $Z_N(\vec{u})$ is of order
  $\frac{\DT}{h^2} h^2$ that is $O(\DT)$ and hence $\|Z_N(\vec{u})\| = O(\DT)$.
  Thus, recalling
  that $\lambda_{\min}(X_N)\geq1$, it holds
  \begin{equation}\label{eq:beta-0}
    \lambda_{\rm min}\left(\frac{A+A^\tr}{2}\right)\geq 1 - \widetilde{C} \DT
  \end{equation}
  for some $\tilde{C}>0$, that contains the infinity norms of the
  first derivatives of $u$ and $D$ and their Lipschitz constants.
  Using Lemma \ref{lem:sigmaA}
  \begin{equation}\label{eq:lambdamin}
    \|A^{-1}\| = \frac{1}{\sigma_{\rm min}(A)} \leq \lambda_{\rm
      min}\left(\frac{A+A^T}{2}\right)^{-1}.
  \end{equation}
  Inequality \eqref{eq:Finvnorm2} now follows combining \eqref{eq:lambdamin} and
  \eqref{eq:beta-0}.

  For the proof of the the estimate in $l^\infty$ norm, we note that
  \[ F^\prime(\vec{u}) = X_N(\vec{u})
  -\frac{\DT}{2h}\tridiag_k\left[u^\prime(\hat\xi_k),O(h),u^\prime(\tilde\xi_k)
  \right]
  \diag_k\left[D^\prime_k\right],
  \]
  where $\hat\xi_k\in[x_{k-1}x_k]$, $\tilde\xi_k\in[x_k,x_{k+1}]$ and
  the constant in the $O(h)$ contains the Lipschitz constant of
  $u^\prime$.
  We split $F^\prime(\vec{u})$ as
  \begin{equation}\label{eq:ZW}
    F^\prime(\vec{u}) =\frac{\DT}{h^2} \left(Z_N-W_N\right),
  \end{equation}
  where
  \begin{align*}
    Z_N &= \diag_k[z_k] =
    \diag_k\left[
      \frac{h^2}{\DT} + D_{k-1/2}+D_{k+1/2}+O(h^2)D^\prime_k
    \right],
    \\
    W_N &= \tridiag_k\left[
      D_{k-1/2}+\frac{h}2u^\prime(\hat\xi_k)D^\prime_{k-1},
      0,
      D_{k+1/2}+\frac{h}2u^\prime(\tilde\xi_k)D^\prime_{k+1}
    \right].
  \end{align*}
  From \eqref{eq:ZW}, we have
  \begin{equation}\label{eq:FZW}
    \left[F^\prime(\vec{u})\right]^{-1}=
  \frac{h^2}{\DT}
  \left(I-Z_N^{-1}W_N\right)^{-1}Z_N^{-1}.
  \end{equation}
  For the factor $Z_N^{-1}$, it holds
    \begin{equation}\label{eq:Zinv}
        \|Z_N^{-1}\|_\infty
  =\max_k
  \frac{1}{|D_{k-1/2}+D_{k+1/2}+\frac{h^2}{\DT}\left(1+O(\DT)D^\prime_k\right)|}
  \leq c \frac{\DT}{h^2},
    \end{equation}
  for $h$ sufficiently small and assuming that $\DT\leq C_\infty h$, recalling
  that $D(\cdot)\geq0$.
  For the remaining factor $(I-Z_N^{-1}W_N)^{-1}$, we note that
  \[ Z_N^{-1}W_N = \tridiag_k
  \left[
    \frac{D_{k-1/2}+\tfrac{h}2u^\prime_kD^\prime_k+O(h^2)}{z_k}
    , \, 0, \,
     \frac{D_{k+1/2}+\tfrac{h}2u^\prime_kD^\prime_k+O(h^2)}{z_k}
   \right]
   \]
   and hence
   \[ \|Z_N^{-1}W_N\|_\infty \leq
   \max_k \frac{D_{k-1/2}+D_{k+1/2}+O(h)}
   {|D_{k-1/2}+D_{k+1/2}+\tfrac{1}{C_\infty}h\left(1+O(\DT)D^\prime_k\right)|}
   \leq\alpha<1,
   \]
   for $C_\infty>0$ sufficiently small.
  Thus the spectral radius of $Z^{-1}W$ is $\rho(Z^{-1}W)<1$
  and we have
  \begin{equation}\label{eq:serie}
    \left(I-Z_N^{-1}W_N\right)^{-1} =
  \sum_{j=0}^\infty \left(Z_N^{-1}W_N\right)^j
  \quad\Rightarrow\quad
  \|\left(I-Z_N^{-1}W_N\right)^{-1}\|_\infty
  \leq \frac{1}{1-\alpha}.
  \end{equation}
  Finally, combining \eqref{eq:serie} and \eqref{eq:Zinv} with
  \eqref{eq:FZW}, the \eqref{eq:Finvnorminf} holds with $C_1=\tfrac{c}{1-\alpha}$.
 \end{proof}

\begin{remark} \label{remark:weaker hp}
The above result, with minor changes, can be proved under weaker assumptions.
Indeed if both $u(\cdot)$ and $D(\cdot)$ are continuously differentiable, then every entry
of $Z_N(\vec{u})$ is of order
\[
O\left(\max\left\{ \frac{\DT}{h}\omega_{u^{\prime}}(h),\frac{\DT}{h}\omega_{D^{\prime}}(h)\|u^{\prime}\|_\infty\right\}\right),
\]
with $\omega_v(\cdot)$ denoting the modulus of continuity of a given function $v$. Therefore with the choice
of $\DT$ proportional to $h$ and setting $\alpha(h)=\max\{\omega_{u^{\prime}}(h),\omega_{D^{\prime}}(h)\}=o(1)$, we find
\[
\lambda_{\rm min}\left(\frac{A+A^\tr}{2}\right)\geq 1 - \widetilde{C}\alpha(h)
\]
and by Lemma \ref{lem:sigmaA} $\|A^{-1}\| \leq 1 +C\alpha(h)$.
Furthermore, if we require that $u$ is only Lipschitz continuous then
the inequality regarding the norm of $A^{-1}$ reads as $\|A^{-1}\|
\leq C$ where $C$ linearly depends on the Lipschitz constant of $u$.
Finally, the same results can be obtained with minor changes, when
using the induced $l^\infty$ norm.
\end{remark}
\begin{remark} \label{remark:real hp} In general, the solution $u$ to
  \eqref{eq:parab} is not smooth, but only piecewise smooth with a
  finite number of cusps. For instance with $D(u)=u^m$ and continuous
  data with piecewise continuous derivative, the derivative of $u$ is
  not defined in a finite number of points in 1D and in a finite
  number of smooth curves in 2D; see \citet{VazquezBOOK}.  The latter
  implies that the related matrices have the same features up to low
  rank correction terms whose cumulative rank is $O(N^{d-1})$ if the
  equation is in $d$ dimensions.
\end{remark}

Since the Crandall-Liggett formula does not induce any restriction on
the the timestep $\DT$ \citep{CL71}, we have only to prove the
convergence of the Newton method. We are interested
in the choice $\DT=Ch$ for a constant $C$ independent of $h$, which
gives a method which is overall first order convergent. This is no
restriction due to the presence of singularities at degenerate points:
higher order methods would be more computationally intensive without
reaching their convergence rate, even if in practice a certain reduction of the
error is expected.

Indeed, concerning the stopping criterion
\( \| \vec{u}^{n+1,s+1} -\vec{u}^{n+1,s} \| \leq \varepsilon \) for
the Newton method, the following observation is of interest. Since the
method is of first order in time $\DT$ can be chosen equal to $h$, it
is sufficient to set $\varepsilon=c\cdot h$ where $c$ is moderately
small constant independent of $h$. In fact, more precision will be
useless in practice and would make the Newton process more expensive
by increasing the iteration count.  The following result is a
classical tool (see \citet{Ortega:book}) for handling the global
convergence of the Newton procedure.

\begin{theorem}[Kantorovich]\label{th:kan}
Consider the Newton method for approximating the zero of a vector
function $F(\vec{u})$, starting from the
initial approximation
$\vec{u}^{(0)}$. Under the assumptions that
\begin{subequations}\label{eq:hp}
\begin{align}
&\|\left[F^{\prime}(\vec{u}^{(0)})\right]^{-1}\| \leq\beta\, ,
\label{eq:hp:beta}\\
&\|\left[F^{\prime}(\vec{u}^{(0)})\right]^{-1}F(\vec{u}^{(0)})\|
\leq\eta\, ,
\label{eq:hp:eta}\\
&\|F^{\prime}(\vec{u})-F^{\prime}(\vec{v})\|\leq\gamma\|\vec{u}-\vec{v}\|\,
, \label{eq:hp:gamma}
\end{align}
\end{subequations}
and that
\begin{equation}\label{eq:12}
    \beta\eta\gamma < \frac12\, ,
\end{equation}
the method is convergent and, in addition,
the stationary point of the iterations lies in the ball with centre
$\vec{u}^{(0)}$ and radius
\[\frac{1-\sqrt{1-2\beta\eta\gamma}}{\beta\gamma}.\]
\end{theorem}

For the choice $\DT=Ch$ we can prove the following result.
\begin{theorem}\label{prop:newton}
  The Newton method for $F(\vec{u})$ defined in \eqref{eq:FnonlinD}
  for computing $\vec{u}^{n}$ is convergent when
  initialised with the solution at the previous timestep (i.e.
  $\vec{u}^{n,0}=\vec{u}^{n-1}$) and for $\DT\leq C h$, for a positive
  constant $C$ independent of $h$.
\end{theorem}
\begin{proof}
  We will make use of the Kantorovich Theorem \ref{th:kan}, so we need
  the estimates \eqref{eq:hp} and to show that \eqref{eq:12} is
  satisfied. We will use the $l^p$ vector norm $\|\vec{u}\|_p^p=\sum
  |v_j|^p$ and the related induced matrix norms. When $p=2$ we find
  the Euclidean vector norm and the induced spectral norm; in general
  they are simply denoted as $\|\cdot\|$.

  Concerning \eqref{eq:hp:beta}, Proposition \ref{prop:Fprimo} and the
  assumption $\DT\leq C_\infty h$ imply
  \begin{equation}\label{eq:beta}
    \beta \leq C_1,
  \end{equation}
  at least for $p=2,\infty$, $C_1=C_1(p)$ and $h$ small enough.

  Regarding \eqref{eq:hp:eta}
  \[
  \left\|\left[F^{\prime}(\vec{u}^{n-1})\right]^{-1}F(\vec{u}^{n-1})\right\|_p
  \leq \beta \left\|F(\vec{u}^{n-1})\right\|_p =
  \beta\left\|\frac{\DT}{h^2}
    L_{D(\vec{u}^{n-1})}\vec{u}^{n-1}\right\|_p =
  \beta\left\|\vec{u}^{n-2}-\vec{u}^{n-1}\right\|_p \leq \beta C_2 \DT
  h^{-1/p}
  \]
  for a constant $C_2=C_2(p)$ independent of $h$.  The first equality
  in the previous calculation follows from \eqref{eq:FnonlinD}, while
  the second one is a consequence of the fact that $\vec{u}^{n-1}$ is
  the stationary point of the Newton iteration for the previous time
  step and thus it satisfies
  \[ \vec{u}^{n-1} + \frac{\DT}{h^2} L_{D(\vec{u}^{n-1})}\vec{u}^{n-1}
  =\vec{u}^{n-2}.\] It follows that
  \begin{equation}\label{eq:eta}
    \eta = C_2\beta\DT h^{-1/p}.
  \end{equation}

  From now on we consider only the $\|\cdot\|_\infty$ norm i.e.
  $p=\infty$, which leads to the most convenient estimate in
  \eqref{eq:eta} and hence to the weakest constraint on the timestep
  $\DT$.

  For the Lipschitz constant of $F^{\prime}$, i.e., for estimating
  \eqref{eq:hp:gamma}, observe that
  $F^{\prime}(\vec{u})-F^{\prime}(\vec{v})$ is a tridiagonal matrix
  with two contributions:
  \begin{equation}\label{eq:g1}
    F^{\prime}(\vec{u})-F^{\prime}(\vec{v})
    = \frac{\DT}{h^2}(L_{D(\vec{u})}-L_{D(\vec{v})})
    + (Y_N(\vec{u})-Y_N(\vec{v})),
  \end{equation}
  with $Y_N(\cdot)$ as in \eqref{Yn}. The first term can be estimated
  as follows:
  \begin{equation}\label{eq:g2}
    \|L_{D(\vec{u})}-L_{D(\vec{v})}\|_\infty \leq
    4\|D^\prime\|_{\infty}\|\vec{u}-\vec{v}\|_\infty.
  \end{equation}
  In order to check that the last inequality is satisfied, one
  observes that the sum of the absolute values of the entries in each
  row of $L_{D(\vec{u})}-L_{D(\vec{v})}$ is smaller than the sum of
  $4$ terms of the form
  \begin{multline*}
    \left|D_{k\pm1/2}(\vec{u})-D_{k\pm1/2}(\vec{v})\right|
    = \left|D\left(\frac{u_{k\pm1}+u_k}2\right) - D\left(\frac{v_{k\pm1}+v_k}2\right)\right| \\
    = \left|D'(\zeta)\right|\frac{|u_{k\pm1}+u_k-v_{k\pm1}-v_k|}2
    \leq\|D'\|_\infty\|\vec{u}-\vec{v}\|_\infty.
  \end{multline*}
  For the second term in \eqref{eq:g1}, we have
  \begin{equation}\label{eq:g3}
    \|Y_N(\vec{u})-Y_N(\vec{v})\|_\infty \leq
    \frac{\DT}{2h^2}\|D^{\prime}\|_\infty\|M\|_\infty,
  \end{equation}
  where
  \[
  M= \tridiag_k\left[
    \begin{aligned}
      &(u_{k-1}-u_k)-(v_{k-1}-v_k)\\
      &(u_{k-1}-2u_k+u_{k+1})-(v_{k-1}-2v_k+v_{k+1})\\
      &(u_{k+1}-u_k) - (v_{k+1}-v_k)
    \end{aligned}
  \right]\]
  and hence
  \begin{equation}\label{eq:g4}
    \|M\|_\infty \leq
    8 \|\vec{u}-\vec{v}\|_\infty.
  \end{equation}
  Replacing equation \eqref{eq:g4} in \eqref{eq:g3} and combining
  \eqref{eq:g3} and \eqref{eq:g2} with \eqref{eq:g1}, we obtain
  \begin{equation}\label{eq:gamma}
    \gamma \leq 8\|D^{\prime}\|_\infty\frac{\DT}{h^2}.
  \end{equation}
  
  Finally, combining equations \eqref{eq:beta}, \eqref{eq:eta}, and
  \eqref{eq:gamma}, Theorem \ref{th:kan} implies that Newton converges
  provided that
  \[ \frac12 \geq  C_1^2C_2 8 \|D^{\prime}\|_\infty\frac{(\DT)^2}{h^2} \geq \beta\eta\gamma,
  \]
  i.e., $\DT\leq C h$, for $h$ sufficiently small and for
  $C=\min\{C_\infty,1/(4C_1\sqrt{C_2\|D^{\prime}\|_\infty\|D^{\prime}\|_\infty}\,)\}$
  (essentially) independent on $h$.
\end{proof}

\section{Algorithms for the resulting linear systems}
\label{sec:LIN}
At each Newton iteration, we need to solve a linear system whose
coefficient matrix is represented by the Jacobian
$F^{\prime}(\vec{u})$ with entries as in \eqref{eq:JacD}. In
principle, the Jacobian is recomputed at each Newton iteration, so we
are interested in efficient iterative methods for solving the related
linear system.

We note in passing that the form of the Jacobian matrix used here is
very similar to the one that is obtained discretizing in space with
$P_1$ conforming finite elements. Thus the methods considered here can
be to some extent generalized to finite elements approximations. In
particular, when considering real $2D$ and $3D$ cases, the structure
of the relevant matrices will depend heavily on the geometry of the
domain, on the triangulation/gridding (often generated automatically),
and on the type of finite elements (higher order or non Lagrangian
etc.). Therefore fast methods that are based on a rigid algebraic
structure (e.g. of Toeplitz type) cannot be adapted because the
structure is lost, in the general framework. However there exists a
kind of information depending only on the continuous operator and
which is inherited virtually unchanged in both finite differences and
finite elements, provided that the grids are quasi-uniform in finite
differences and the angles are not degenerating in finite
elements. Such information consists in the locally Toeplitz structure
(see \citet{Serra06:glt,Tilli98}) and in the related spectral features
(conditioning, subspaces related to small eigenvalues etc.). We remind
that these spectral features are conveniently used when defining ad
hoc preconditioned Krylov methods or multigrid algorithms, working
uniformly well in one or more dimensions.

In order to choose appropriate iterative methods for solving the
jacobian linear system, we first analyse the spectral properties of
the matrix $F^\prime(\vec{u})$. This will lead us to consider
preconditioned Krylov methods, multigrid and their combinations.

\subsection{Spectral analysis for the resulting matrix-sequences}
\label{ssec:SPECTRAL}

We start by introducing the notion of spectral distribution for a
matrix sequence. Then we will briefly report a concise analysis of
some delicate spectral features of the matrices involved in the
definition of the Jacobian. Since the emphasis of this work relies
in the computational aspects, we will not report all possible
details, nuances, and generalisations of the spectral analysis.

\begin{definition}\label{def-distribution}
Let $\mathcal C_0({\mathbb R}^+_0)$ be the set of continuous
functions with bounded support defined over the nonnegative real
numbers, $d$ a positive integer, and $\theta$ a complex-valued
measurable function defined on a set $G\subset\mathbb R^d$ of finite
and positive Lebesgue measure $\mu(G)$. Here $G$ will be often equal
to $(-\pi,\pi)^d$ so that $e^{\imm\overline{G}}={\mathbb T}^d$ with
$\imm^2 = -1$ and
${\mathbb T}$ denoting the complex unit circle.
A matrix sequence $\{A_N\}$ is said to be {\em
distributed $($in the sense of the eigenvalues$)$ as the pair
$(\theta,G)$,} or to {\em have the eigenvalue distribution function
$\theta$} ($\{A_N\}\sim_{\lambda}(\theta,G)$), if, $\forall F\in
\mathcal C_0({\mathbb C})$, the following limit relation holds
\begin{equation}\label{distribution:eig}
\lim_{N\rightarrow \infty}{\frac 1 {N}} \sum_{j=1}^{N}
F\left(\lambda_j(A_N)\right)=\frac1{\mu(G)}\,\int_G F(\theta(t))\,
dt,\qquad t=(t_{1},\ldots,t_{d}).
\end{equation}
Furthermore, a matrix sequence
$\{A_N\}$ is said to be {\em distributed $($in the sense of the
singular values$)$ as the pair $(\theta,G)$,} or to {\em have the
distribution function $\theta$} ($\{A_N\}\sim_{\sigma}(\theta,G)$),
if, $\forall F\in \mathcal C_0({\mathbb R}^+_0)$, the following
limit relation holds
\begin{equation}\label{distribution:sv-eig}
\lim_{N\rightarrow \infty}{\frac 1 {N}} \sum_{j=1}^{N}
F\left(\sigma_j(A_N)\right)=\frac1{\mu(G)}\,\int_G F(|\theta(t)|)\,
dt,\qquad t=(t_{1},\ldots,t_{d}).
\end{equation}
\end{definition}

Along with the distribution in the sense of singular
values/eigenvalues (weak*-convergence), for the practical
convergence analysis of iterative solvers we are also interested in
a further asymptotic property called here the {\em clustering}.
\begin{definition}\label{def-cluster}
A matrix sequence $\{A_N\}$ is {\em strongly clustered at $s \in
\mathbb C$} (in the eigenvalue sense), if for any $\varepsilon>0$
the number of the eigenvalues of $A_N$ off the disk
\[
D(s,\varepsilon):=\{z:|z-s|<\varepsilon\}
\]
can be bounded by a pure constant $q_\varepsilon$ possibly depending on
$\varepsilon$, but not on $n$. In other words
\[
q_\varepsilon(n,s):=\#\{\lambda_j(A_N): \lambda_j\notin
D(s,\varepsilon)\}=O(1), \quad n\to\infty.
\]
If every $A_N$ has only real eigenvalues (at least for all $n$ large
enough), then $s$ is real and the disk $D(s,\varepsilon)$ reduces to
the interval $(s-\varepsilon,s+\varepsilon)$. Furthermore, $\{A_N\}$ is
{\em  strongly clustered at a nonempty closed set $S \subset \mathbb
C$} (in the eigenvalue sense) if for any $\varepsilon>0$
\begin{equation} \label{2.2}
q_\varepsilon(n,S):=\#\{\lambda_j(A_N): \lambda_j\not\in
D(S,\varepsilon):=\cup_{s\in S} D(s,\varepsilon)\}=O(1), \quad n\to\infty,
\end{equation}
$D(S,\varepsilon)$ is the $\varepsilon$-neighbourhood of $S$, and if every
$A_N$ has only real eigenvalues, then $S$ has to be a nonempty
closed subset of $\mathbb R$. Finally, the term ``strongly'' is
replaced by ``weakly'', if
\[
q_\varepsilon(n,s)=o(n), \qquad \bigl(q_\varepsilon(n,S)=o(n)\bigr), \quad
n\to\infty,
\]
in the case of a point $s$ (a closed set $S$), respectively. The
extension of the notion in the singular value sense is trivial and
is not reported in detail.
\end{definition}
\begin{remark}
  It is clear that $\{A_N\}\sim_{\lambda} (\theta,G)$
  ($\{A_N\}\sim_{\sigma} (\theta,G)$) with $\theta\equiv s$ a constant
  function is equivalent to $\{A_N\}$ being weakly clustered in the
  eigenvalues sense at $s \in \mathbb C$ (in the singular value sense
  at $s \in {\mathbb R}^+_0$).
\end{remark}

Now we briefly use the above concepts in our specific setting. Given
the linear restriction on $\DT$ imposed by the convergence of the
Newton method (Theorem \ref{prop:newton}), we are interested in the
choice $\DT=Ch$ for $C>0$ independent of $h$. However, for
notational simplicity, here we assume $\DT=h$ and note that
analogous results hold for $C>0$.

Taking into account $\DT=h$ and the re-scaling
$A_N=h F^\prime(\vec{u})$, we consider the sequence
$\{A_N\}$ such that
\begin{eqnarray*}
A_N & = & - L_{D(\vec{u})} + R_N(\vec{u}) \\
R_N(\vec{u}) & = & h I_N - \frac12 T_N(\vec{u})\diag_k(D^\prime_k)
\end{eqnarray*}
with $T_N$ defined as in \eqref{Tn}.

We have the following results, which are of crucial interest in the
choice, in the design, and in the analysis of efficient solvers for
the involved linear systems.

\begin{remark}
  \label{rem:condizionamento}
  The conditioning in spectral norm of $A_N$ is of order $N$: this is
  implied directly by Proposition \ref{prop:Fprimo}.
  More in detail, by
  using the Bendixson Theorem \citep[see][Theorem 3.6.1]{Sto02} the eigenvalues of $A_N$ are localised
  in a rectangle having real part in $[c h, C]$ and imaginary part in
  $[-d h, d h]$ for some positive constants $c$, $d$, $C$ independent of
  $N$. This statement is again implied by the analysis provided in
  Proposition \ref{prop:Fprimo} for the real part, while for the imaginary
  part we note that $(A-A^T)/2 = -\frac{1}{4}
  \tridiag_k[(D^{\prime}_{k-1}+D^{\prime}_k)(u_{k-1}-u_k), \; 0, \;
  (D^{\prime}_{k+1}+D^{\prime}_k)(u_{k+1}-u_k)]$.
\end{remark}

\begin{remark}
  \label{rem:laplaciano}
  $\{A_N\}\sim_{\lambda, \sigma}(\theta,G)$ with
  $\theta(x,s)=D(u(x))(2-2\cos(s))$, $G=[a,b]\times [0,2\pi]$
  (distribution of the zero order main term). The distribution of
  $\{L_{D(\vec{u})}\}$ is already known \citep[see][]{Tilli98}, if we
  assume that $\vec{u}$ is a sampling of a given function over a
  uniform grid. In our case the entries of $\vec{u}$ represent an
  approximation in infinity norm of the true solution, the latter
  being implied by the convergence of the method, and therefore by
  standard perturbation arguments we deduce
  $\{L_{D(\vec{u})}\}\sim_{\lambda, \sigma}(-\theta,G)$ with $\theta$
  and $G$ as above.  Moreover the trace norm (sum of all singular
  values i.e. Schatten $p$ norm with $p=1$; see \citet{Bhatia:book})
  of the remaining part $R_N(\vec{u})$ is bounded by a pure constant
  $C$ independent of $N$, when assuming that $D'$ is bounded and $u$
  is at least Lipschitz continuous. The latter implies that the
  distribution of $\{A_N\}$ is decided only by the symmetric part that
  is, essentially, $\{L_{D(\vec{u})}\}$ \citep[see][Theorem
  3.4]{GS07:jacobi-nonsymm}. Moreover any real interval containing the
  spectrum $\{L_{D(\vec{u})}\}$ is also a strong eigenvalue clustering
  set for $\{A_N\}$ \citep[see][Corollary 3.3 and Theorem
  3.5]{GS07:jacobi-nonsymm}.

  Concerning the negligible term, we have that
  $\{R_N(\vec{u})\}\sim_{\lambda, \sigma}(0,G)$ and
  $\{R_N(\vec{u})/h\}\sim_{\lambda, \sigma}(\psi,G)$ with
  $\psi(x,s)=1-D'(u(x))2\imm \sin(s))$, $G=[a,b]\times [0,2\pi]$
  (distribution of the first order term).
\end{remark}

\begin{remark}
  \label{rem:precondizionamento}
  Setting
  \[P_N=- L_{D(\vec{u})} + h I_N,\] we have
  $\{P_N^{-1}A_N\}\sim_{\lambda, \sigma}(1,G)$ (equivalent, as already
  observed, to a weak eigenvalue/singular value clustering): it
  follows from the property of algebra of the Generalized Locally
  Toeplitz (GLT) sequences \citep[see][]{Serra06:glt}.

  In fact the preconditioned sequence $\{P_N^{-1}A_N\}$ is also
  strongly clustered at $1$ both in the eigenvalue and singular value
  sense: we remark that the strong clustering property can be
  recovered via local domain analysis, by employing the same tools and
  the same procedure as in \citet[Theorem 3.7]{BGST05}; see also
  Section 3.1 and the conclusion section in \citet{BS07} and
  references therein.

  More in detail, by the Bendixson Theorem the eigenvalues of
  $P_N^{-1}A_N$ are localised in a rectangle having real part in
  $[1-c_1 h, \, 1+c_2h]$ and imaginary part in $[-d, \, d]$ for some
  positive constants $c_1$, $c_2$, $d$ independent of $N$.  This
  statement follows by noting that the eigenvalues of $P_N^{-1}A_N$
  belong to the field of value of
  $P_N^{-1/2}A_NP_N^{-1/2}$. Considering $\alpha =
  \vec{x}^HP_N^{-1/2}A_NP_N^{-1/2}\vec{x}$, for all $\vec{x} \in
  \C^n$, $\|\vec{x}\|=1$, it holds that the real part of $\alpha$ is
  $\vec{x}^HP_N^{-1/2}(A_N+A_N^T)P_N^{-1/2}\vec{x}/2$ which belongs to
  $[1-c_1 h, \, 1+c_2h]$ by the analysis provided in Proposition
  \ref{prop:Fprimo}. A similar analysis stands for the imaginary part
  of $\alpha$ similarly to Remark~\ref{rem:condizionamento}.
\end{remark}

Remark \ref{rem:precondizionamento} is very important in practice,
since it is crucial for deducing that the number of iterations of
preconditioned GMRES is bounded by a constant depending on the precision,
but not on the mesh that is on $h$ (optimality of the method). This
will be discussed in the next section.

\subsection{Iterative methods for the linear system}
\label{ssec:MULTI}
In this section we consider some iterative methods for solving the
linear system at each Newton step and study their convergence
properties on the matrix sequence $\{A_N\}$. A classical reference for
the results quoted below is \citet{Saad:book}.

\paragraph{GMRES}
We first consider the GMRES algorithm, since the antisymmetric part of
$A_N$ is negligible but not zero.

Assume that $A_N$ is diagonalisable and let $A_N=W\Lambda{W}^{-1}$,
where $\Lambda=\diag_k(\lambda_k)$ is the diagonal matrix of the
eigenvalues.  Define
\[ \epsilon^{(m)} = \min_{p\in\mathbb{P}_m: p(0)=1}
\max_{k=1,\ldots,N} |p(\lambda_k)|.
\]
Denoting with $r^{(m)}$ the residual at the $m^{\text{th}}$ step of
GMRES, it is a classical result that
\[ \| r^{(m)}\|_2 \leq \kappa_2(W) \epsilon^{(m)} \|r^{(0)}\|_2.
\]

Thanks to Remark \ref{rem:condizionamento},
$\kappa_2(W)\approx1$. Thus the GMRES convergence is determined by the
factor \(\epsilon^{(m)}\).  Thanks to Remark
\ref{rem:condizionamento}, it is possible to construct an ellipse
properly containing the spectrum of $A_N$ and avoiding the complex
$0$, so that when one applies GMRES to the matrix $A_N$, it holds that
\begin{equation}\label{eq:GMRESi}
\epsilon^{(m)} \leq \left(1-C\sqrt{h}\right)^m
\end{equation}
for a positive constant $C$ that is independent of the problem size
$N$.

Similarly, using $P_N$ as preconditioner, Remark
\ref{rem:precondizionamento} implies that
\begin{equation}\label{eq:GMRESii}
\epsilon^{(m)} \leq \widetilde{C}^m
\end{equation}
for some $\widetilde{C}\in(0,1)$, independent of the
problem size $N$. Even if the solution $u$ is not enough
regular to assure that the spectrum of $P_N^{-1}A_N$ belongs to
$[1-c_1 h, \, 1+c_2h] \times \imm [-d, \, d]$, the strong cluster at $1$
leads in practice to the super-linear convergence.

\paragraph{Conjugate gradient (CG)}
Let \(S_N=(A_N+A^\tr_N)/2\) be the symmetric part of $A_N$ and define
$\|\vec{x}\|_{S_N}=\|(S_N)^{1/2}\vec{x} \|_2$.  Denoting
$\kappa_2(S_N)=\|S_N\|_2\|S_N^{-1}\|_2$, we recall the following
classical result about the convergence of the CG:
\begin{equation}\label{eq:CG}
  \|\vec{x}_m-\vec{x}_*\|_{S_N}
  \leq 2 \left( \frac{\sqrt{\kappa_2(S_N)}-1}{\sqrt{\kappa_2(S_N)}+1} \right)^m
  \|\vec{x}_0-\vec{x}_*\|_{S_N},
\end{equation}
where $\vec{x}_m$ is the approximate solution obtained at the
$m^{\text{th}}$ step of the CG algorithm and $\vec{x}_*$ the exact
solution.

Thus, combining \eqref{eq:CG} with Remark \ref{rem:condizionamento},
we expect the CG algorithm to converge in $O(\sqrt{N})$ iterations
when applied to $S_N$.  On the other hand, using $P_N$ as
preconditioner, \eqref{eq:CG} together with Remark
\ref{rem:precondizionamento} imply that CG converges in a
constant number of iterations, independently on the size $N$ of the
problem.

Finally, according to Remarks \ref{rem:condizionamento} and
\ref{rem:laplaciano}, the antisymmetric part of $A_N$ is
negligible. Thus in practice one may apply the CG algorithm to the
matrix $A_N$, expecting a convergence behaviour similar to that for
$S_N$, in both the unpreconditioned and preconditioned cases.

\paragraph{Multigrid method (MGM)}

From Remark \ref{rem:laplaciano} we have that $A_N$ has the same
spectral behaviour of $-L_{D(\vec{u})}$. Hence, if an iterative method
is effective for $L_{D(\vec{u})}$ and robust, it should be effective
also for $A_N$. This is the case of MGM largely used with
elliptic PDEs \citep{Trottenberg:book}.

MGM has essentially two degrees of indetermination: the choice of the
grid transfer operators and the choice of the smoother (pre- and
post-smoother, if necessary).
In particular, let $P_{i+1}^i$ be the prolongation operator from a coarse grid
$i+1$ to a finer grid $i$. We consider a Galerkin strategy: the
restriction operator is $(P_{i+1}^i)^\tr$ and the coefficient matrix
of the coarse problem is $A_{i+1}=(P_{i+1}^i)^\tr A_i P_{i+1}^i$,
where $A_i$ is the coefficient matrix on the $i^{\text{th}}$ grid.

For the prolongation we consider the
classical linear interpolation. We note that it is not necessary to
resort to more sophisticated grid transfer operators since $A_N$ is
spectrally distributed as $-L_{D(\vec{u})}$. The restriction is the
full-weight since, according to the Galerkin approach, it is the
transpose of the linear interpolation. Concerning the smoother damped
Jacobi, damped Gauss-Seidel and red-black Gauss-Seidel are considered.

\begin{remark}
  The robustness of MGM could be improved in several way. A
  possibility is to use as post-smoother a damped method that reduces
  the error in the middle frequencies whose could be not well dealt
  with the pre-smoother and the coarse grid correction. This is called
  as ``intermediate iteration'' in the multi-iterative methods
  \citep{Serra93:multi}. Another degree of freedom is the number of
  smoothing iterations depending on the grid $i$. Indeed in
  \citet{ST04} it is shown that a polynomial growth with $i$ does not
  affect the global cost, that remains linear for banded structures,
  only changing the constants involved in the big $O$.
\end{remark}

In our present setting is not necessary to resort to the strategies
described in the previous remark. In fact the method that achieves the
smallest theoretical cost and that minimises the CPU times, for
reaching the solution with a preassigned accuracy $\epsilon$, is the
simplest V-cycle with only one step of damped Jacobi as
pre-smoother. The reason of the observed behaviour relies in the
spectral features of our linear algebra problem: indeed, $A_N$ can be
viewed, after re-scaling, as a regularised weighted Laplacian since in
the coefficient matrix one adds $h$ times the identity (see the
previous subsection). In this way the conditioning is not growing as
$N^2$ as in the standard Laplacian but grows only linearly with $N$
(see Remark \ref{rem:condizionamento}).

Therefore the basic V-cycle, with one single step of damped Jacobi as
pre-smoother, is already optimal for $A_N$, i.e. the number of
iterations is independent of the system size \citep{Trottenberg:book}.
Moreover, as we will see in the numerical tests of section
\ref{ssec:linear-sys}, the number of iterations for reaching a given
accuracy is already very moderate. Therefore the additional cost per
iteration, that should be paid for increasing the number of smoothing
steps and for the use of a post-smoother, can not be compensated by a
remarkable reduction of the iteration count.

Finally, we stress that a robust and effective strategy is to use a
multigrid iteration as preconditioner for GMRES as confirmed in the
numerical experiments. In fact we showed that $P_N$ is an optimal
preconditioner for GMRES and the MGM is an optimal solver for a linear
system with matrix $P_N$.

\section{Numerical tests}
\label{sec:NUM}

In this section we consider as a test case the porous medium equation
written in the form
\begin{equation}\label{eq:PM}
\pder{u}{t} = \pder{}{x}\left(mu^{m-1}\pder{u}{x}\right)
\end{equation}
with homogeneous Dirichlet boundary conditions. Here $m\geq1$, with
$m=1$ corresponding to the heat equation. In particular we
consider the exact self-similar solution
\begin{equation}\label{eq:BB}
  u(t,x) =t^{-\alpha}
  \left[1 - \alpha\tfrac{m-1}{2m}
    \left(|x|t^{-\alpha}\right)^2
  \right]_+^{\frac1{m-1}} ,
  \qquad
  \alpha=\tfrac1{m+1}
\end{equation}
due to Barenblatt and Pattle \citet{VazquezBOOK}. (The subscript $+$
denotes the positive part). The experiments are carried out in
Matlab~7.5.

\subsection{Convergence of the global method and of Newton's method}
\label{ssec:conv-pde-newton}

First we check the convergence of the method. We perform test for
$m$ ranging from $2$ to $5$, observing no appreciable difference in
the convergence properties of the algorithm. In all tests we choose
$\DT=h$.

\begin{figure}
  \centering{
    \includegraphics[width=.4\textwidth]{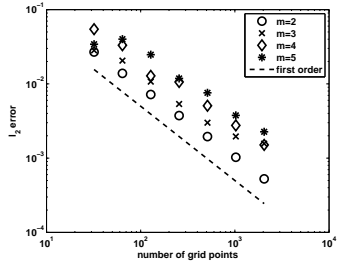}
  }
  \caption{$l_2$ error at final time for $N=32,64,\ldots,2048$, $m=2$,
    final time $t=20/32$, $\DT=h$.}
  \label{fig:convergence}
\end{figure}

Figure \ref{fig:convergence} plots the $l_2$ errors between the
numerical solution at time $t=20/32$ and the exact solution
\eqref{eq:BB} and  shows that the method is first order convergent,
as expected for this choice of time stepping procedure and also due
to the presence of the singularity in the first derivative of the
exact solution. The dashed line is a reference slope for first order
schemes. We observe that the convergence is not significantly
affected by the parameter $m$.

\begin{figure}
  \centering{
    \begin{tabular}{cc}
      \includegraphics[width=.4\textwidth]{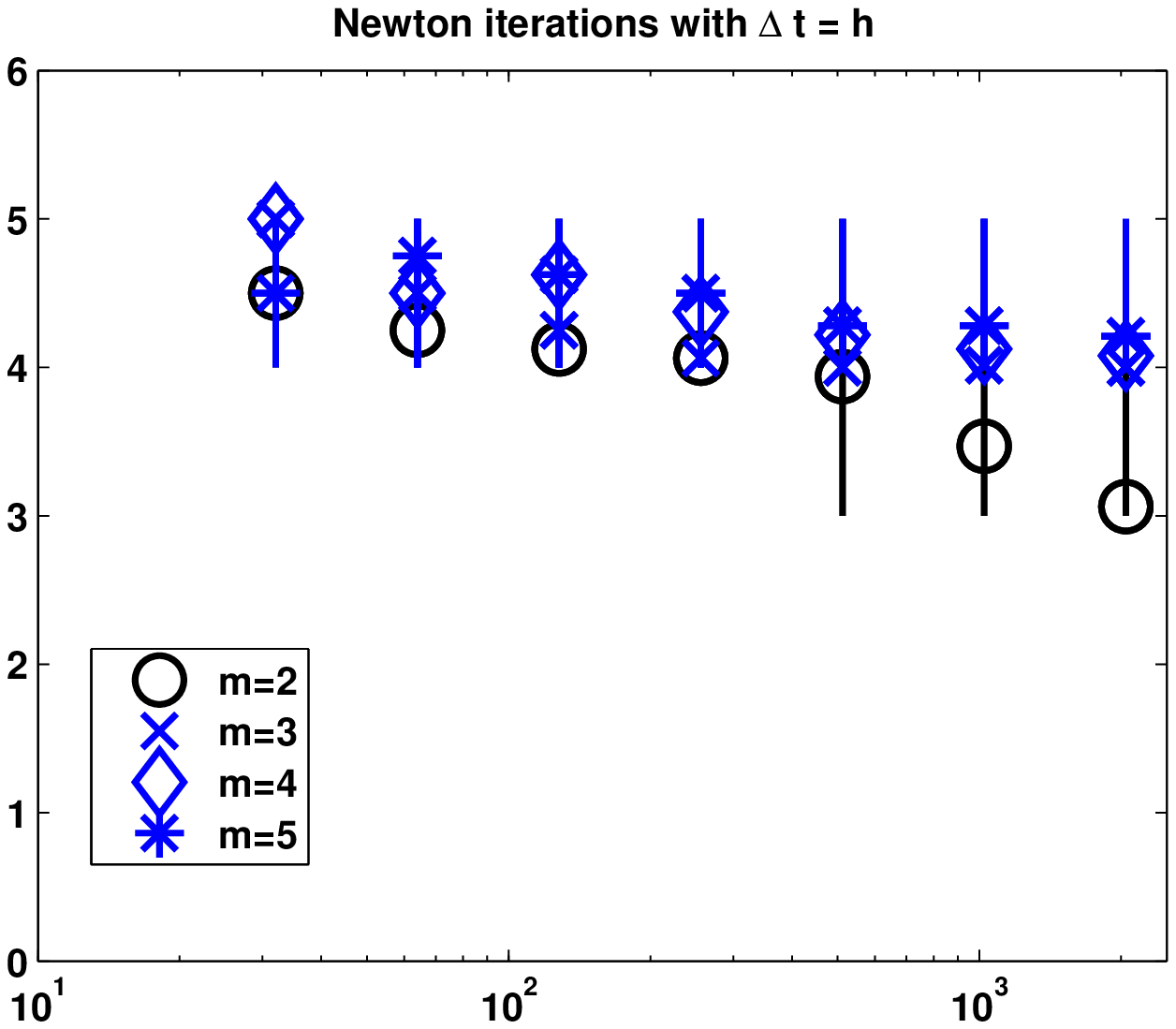}
      &
      \includegraphics[width=.4\textwidth]{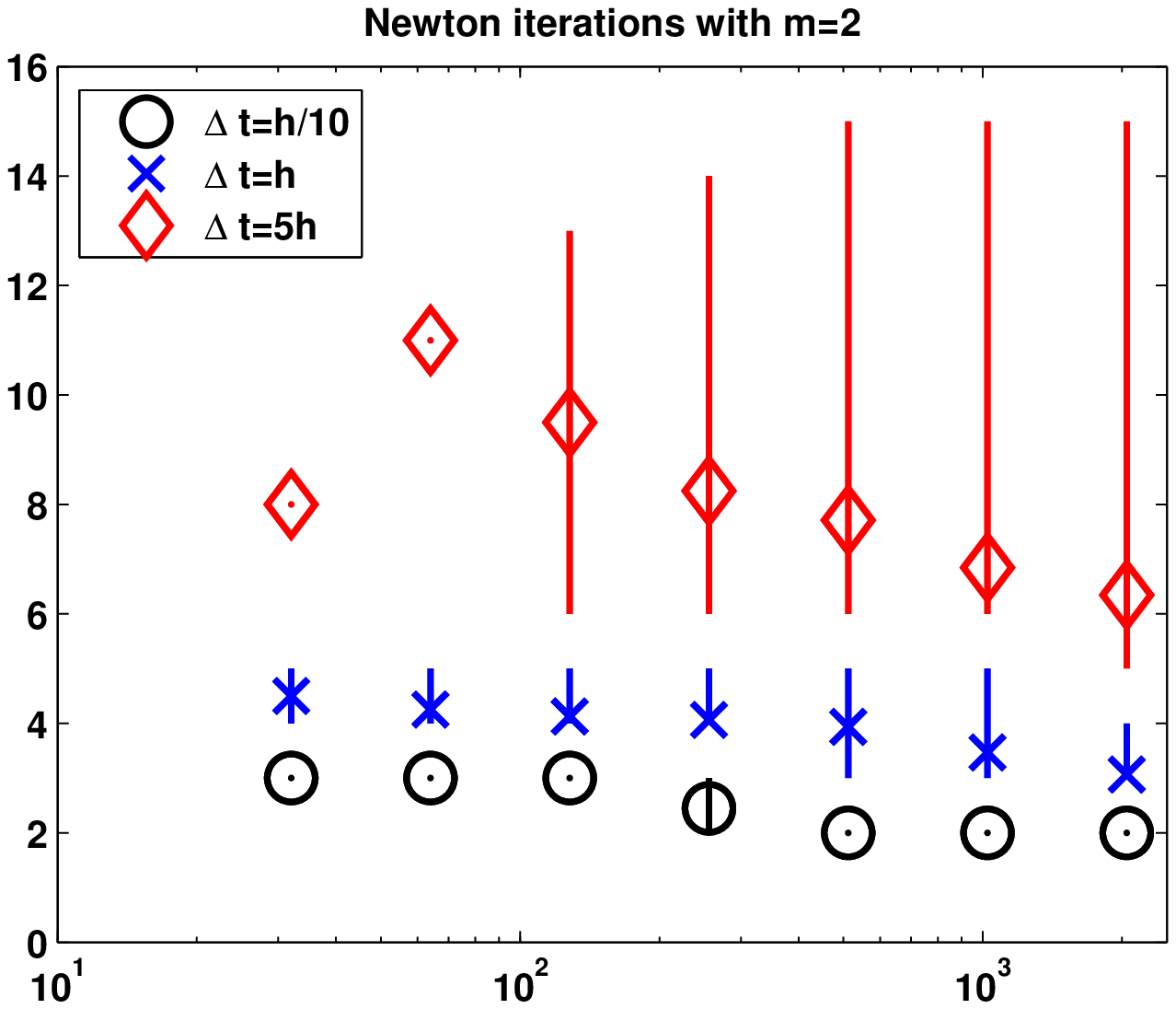}
      \\
      {\bf (a)} & {\bf (b)}
    \end{tabular}
  }
  \caption{Average, minimum and maximum number of Newton iterations
    performed during the integration until final time. In (a) $\DT$
    was kept fixed and $m$ varied, in (b) $m=2$ was kept fixed and
    $\DT$ varied.}
\label{fig:newtoniter}
\end{figure}

Figure \ref{fig:newtoniter} plots the number of Newton iterations
employed by the algorithm during the integration from $t=0$ to
$t=20/32$. We plot the average (circles), minimum and maximum
(solid lines) number of Newton iterations per timestep.  Taking
$\DT=h$ (Figure \ref{fig:newtoniter}a), we observe that the number of
Newton iterations slowly decreases when $N$ increases and that, for
any given $N$ it increases only very moderately when $m$ increases.
In the case $m=2$ we also tried to vary the step size from $\DT=h/10$
to $\DT=5h$. The results are reported in Figure
\ref{fig:newtoniter}b, showing that the number of Newton iterations
grows when taking larger $\DT$ in \eqref{CLformula}. The larger
variability (for fixed $N$) and the irregular behaviour of the mean
value when increasing $N$ in the case $\DT=5h$ preludes to the loss of
convergence that we observe if $\DT$ is taken even larger.

\begin{figure}
  \centering{
    \includegraphics[width=.6\textwidth]{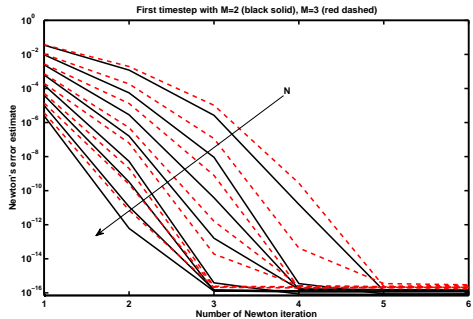}
  }
  \caption{History of the convergence of the Newton iterations during
    the first timestep. Black solid lines correspond to $m=2$ and red
    dashed ones to $m=3$. We show the results for $N$ ranging from
    $32$ to $4096$, with $\DT=h$: the behaviour under grid refinement
    is indicated by the thin arrow.}
  \label{fig:primopasso}
\end{figure}

Next we verify the convergence of the Newton's method. In Figure
\ref{fig:primopasso} we plot the Newton's error estimate
$\|\vec{u}^{1,k+1}-\vec{u}^{1,k}\|/\|\vec{u}^{1,k}\|$ obtained when
computing the first timestep $\vec{u}^1$. We compare different number
of grid points ($N=32,64,\ldots,4096$) as indicated by the thin arrow
and two values for the exponent $m$ appearing in \eqref{eq:PM}.

We emphasise that as prescribed in Proposition \ref{prop:newton} the
choice of $\DT=h$ is acceptable for the convergence both of the
global numerical scheme and for the convergence of the Newton
procedure.

\subsection{Solution of the linear system}\label{ssec:linear-sys}

This section is devoted to computational proposals for the solution of
a linear system where the coefficient matrix is the Jacobian in
\eqref{Jf}, which is required at every step of the Newton procedure.
 For all the tests, we set $m=2$, final time
$t=20/32$, $\DT=h$, and we let $N$ be equal to $32,64,\ldots,1024$
for checking the optimality of the proposed best solvers.

As already stressed in Remark \ref{remark:Yn}, the matrix is
(weakly) non-symmetric so we start by considering the use of
preconditioned GMRES (PGMRES).

\subsubsection{GMRES}
\label{sssec:num:GMRES}
\begin{figure}
  \centering{
    \begin{tabular}{cc}
      \includegraphics[width=.4\textwidth]{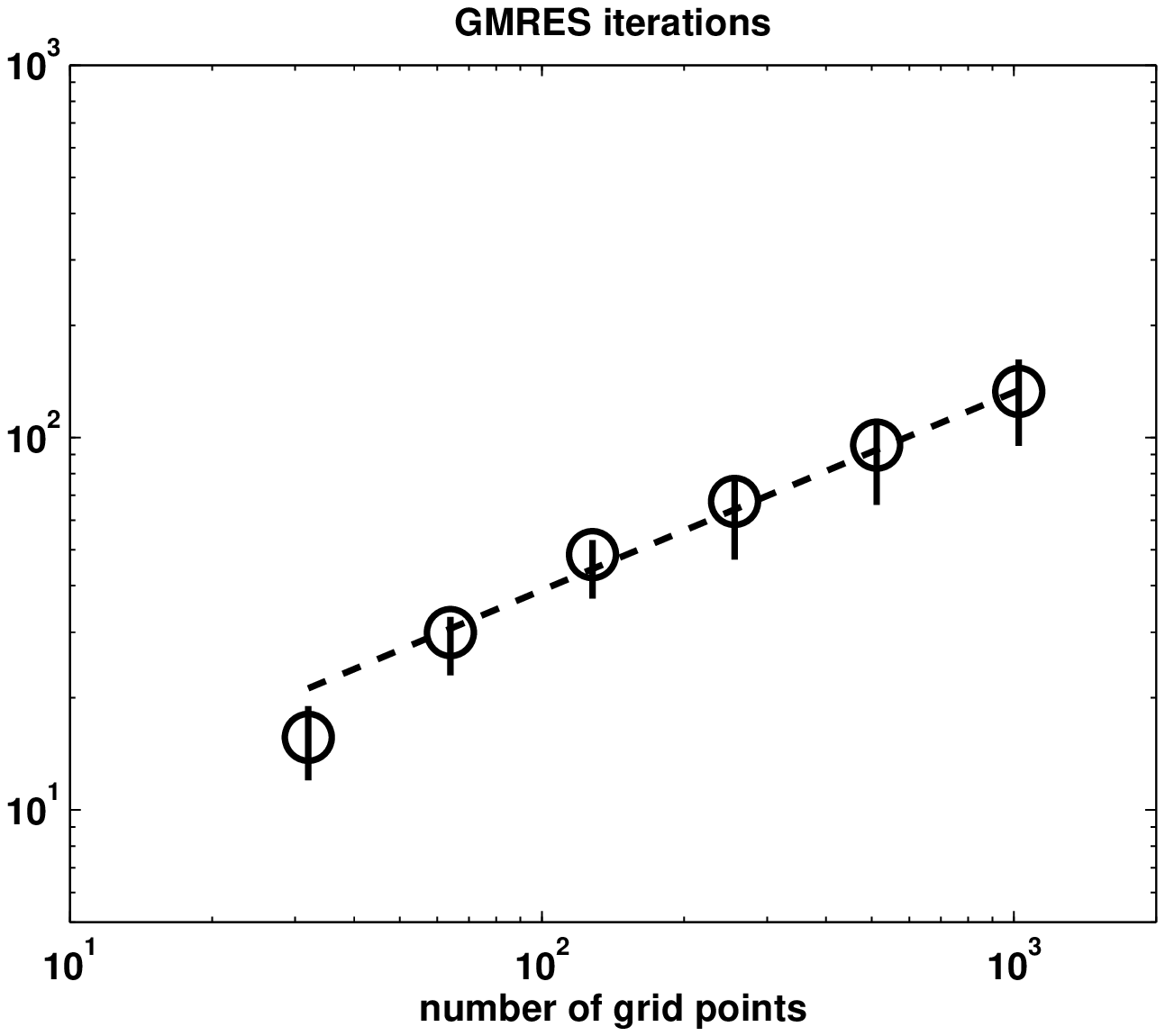}
      &
      \includegraphics[width=.4\textwidth]{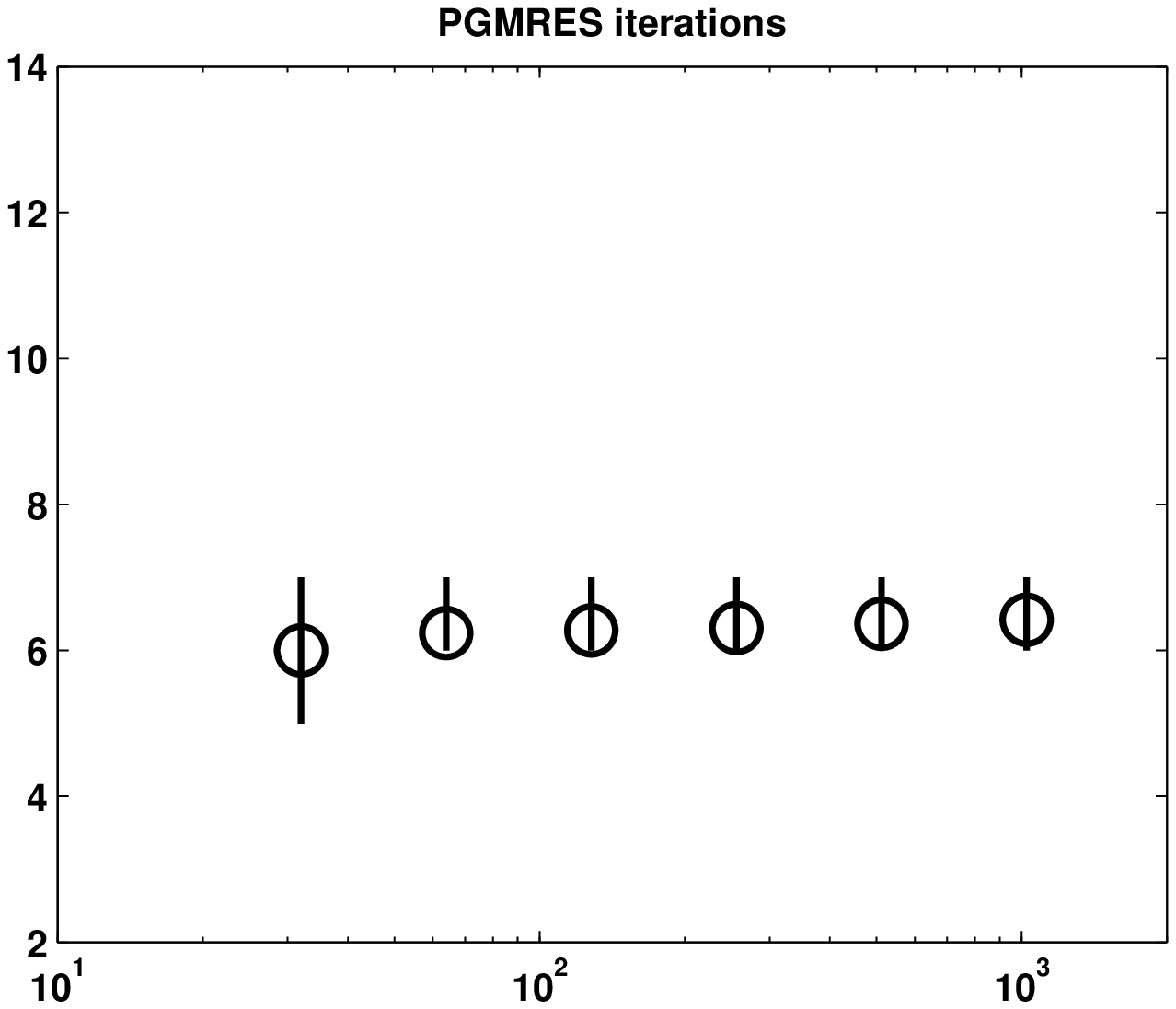}
      \\
      {\bf (a)} & {\bf (b)}
    \end{tabular}
  }
  \caption{Average, minimum and maximum number of GMRES iterations (a)
    and preconditioned GMRES iterations (b) performed during the
    integration until final time. The dashed line in panel (a) is the
    least square fit.
  } \label{fig:GMRES}
\end{figure}

In Figure \ref{fig:GMRES}a we plot the average (circles), minimum and
maximum (vertical lines) number of GMRES iterations performed during
the integration until  final time, at different spatial
resolutions. A least square fit (dashed line) shows that the number of iterations
grows as $N^{0.5320}$. This fact is in complete accordance with the
analysis of Subsection \ref{ssec:MULTI} and in particular with
equation \eqref{eq:GMRESi}.


In Remark \ref{rem:laplaciano} we proved that $Y_N(\vec{u})$ is
negligible with respect to the symmetric positive definite term
$X_N(\vec{u})$. Accordingly, in Remark \ref{rem:precondizionamento}
the use of $X_N(\vec{u})$ as preconditioner for $F^\prime(\vec{u})$
was analysed and it was shown to provide a strong spectral clustering
of the preconditioned matrix at $1$ and therefore we expect a number
of iterations not depending on the size $N$ of the matrix as in
\eqref{eq:GMRESii}: this fact is observed in practice and indeed the
iteration count of the PGMRES is almost constant, with
average value equal to 6 iterations (see Figure \ref{fig:GMRES}b).

At this point we are left with the problem of solving efficiently a
generic linear system with coefficient matrix $X_N(\vec{u})$, which is
a regularised version of a weighted Laplacian (i.e., by
re-scaling, it is a shift of $-L_{D(\vec{u})}$ by $h^2/\DT$ times the
identity). A standard V-cycle is thus optimally convergent since
$X_N(\vec{u})$ is slightly better conditioned than a standard
Laplacian.

\subsubsection{CG}
\label{sssec:num:CG}
\begin{figure}
  \centering{
    \begin{tabular}{cc}
      \includegraphics[width=.4\textwidth]{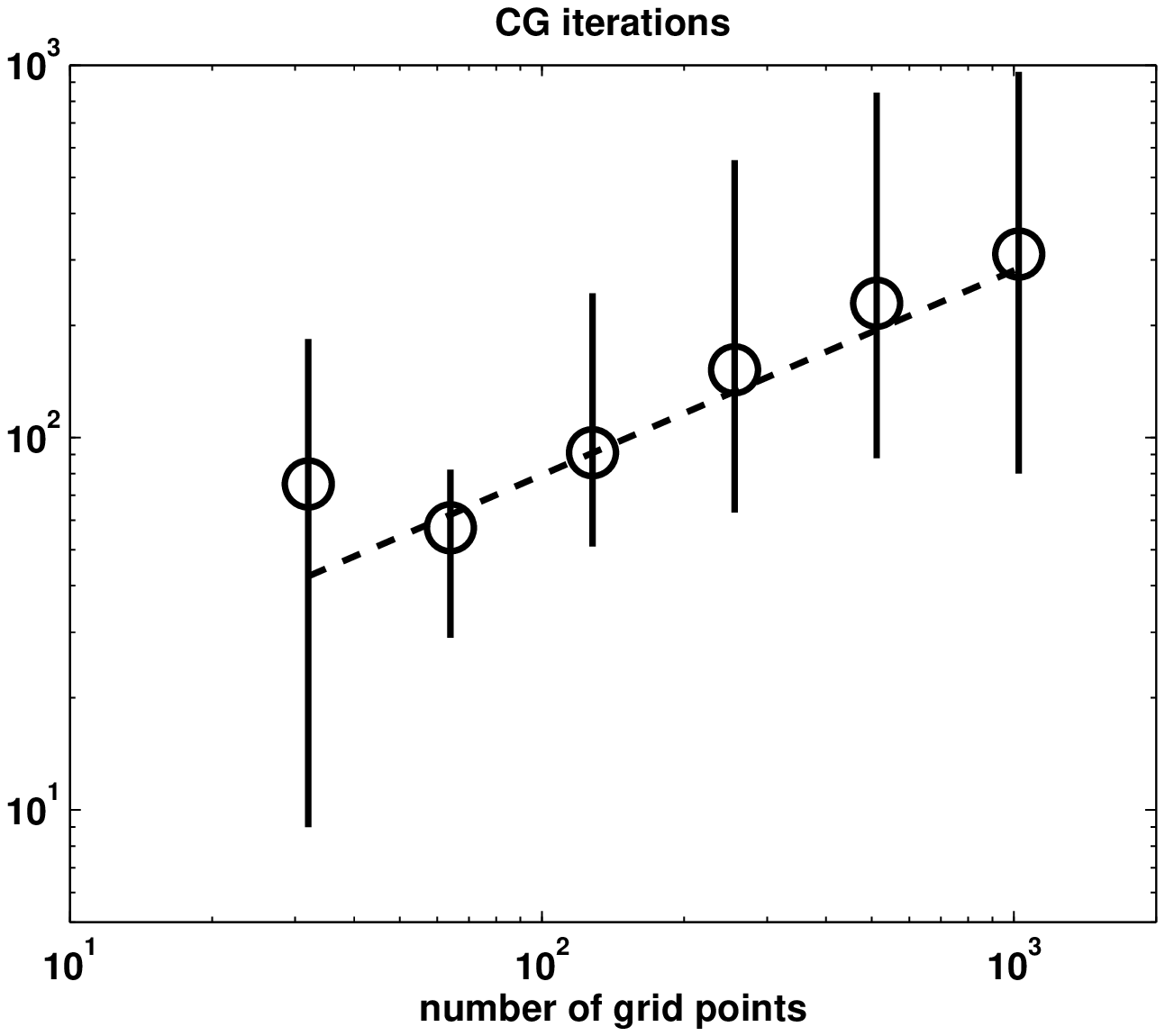}
      &
      \includegraphics[width=.4\textwidth]{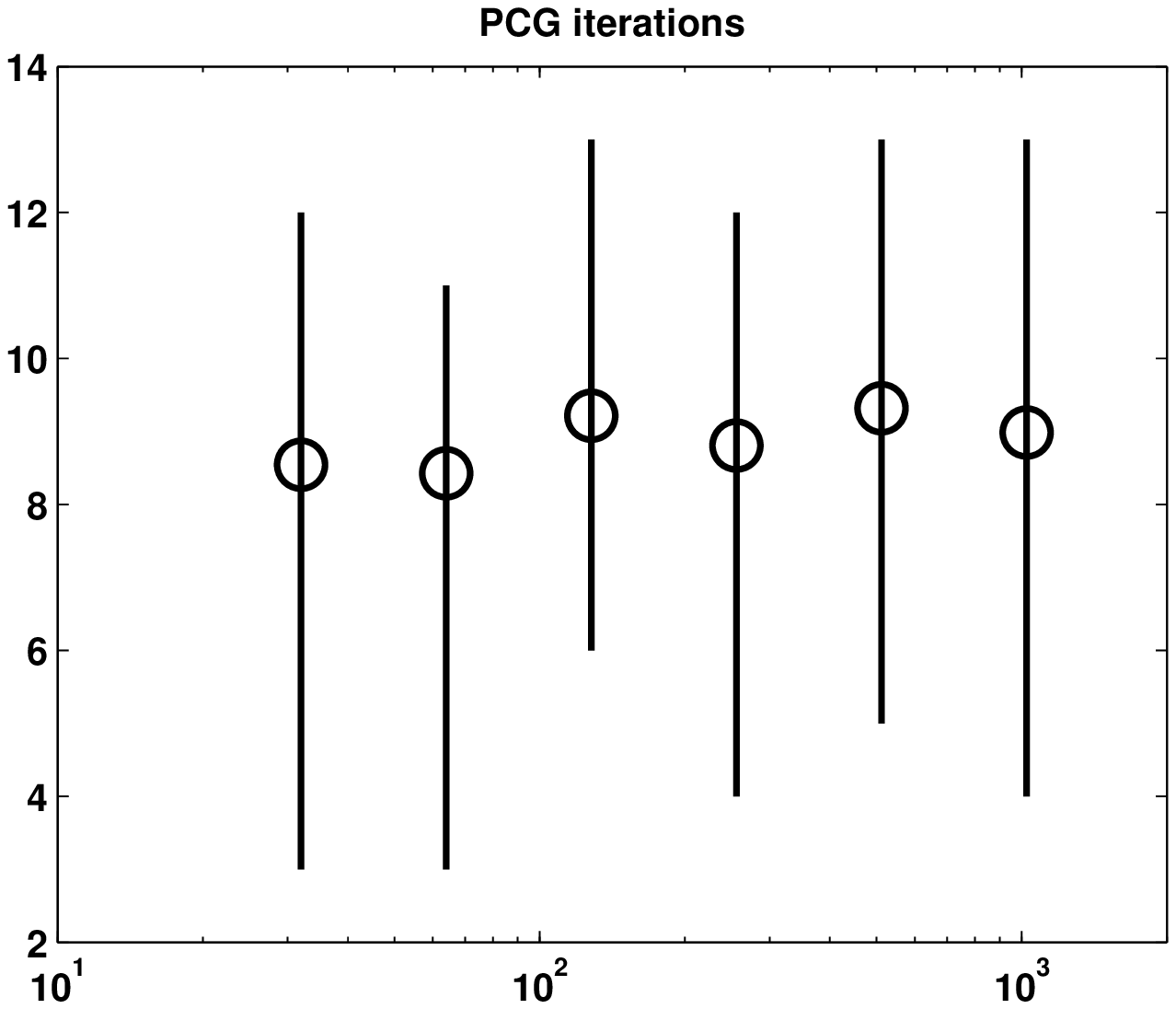}
      \\
      {\bf (a)} & {\bf (b)}
    \end{tabular}
  }
  \caption{Average, minimum and maximum number of CG iterations (a)
    and PCG iterations (b), performed during the
    integration until final time. The dashed line in panel (a) is the
    least square fit.
  } \label{fig:CG}
\end{figure}

Since the non-symmetric part of $F^\prime(\vec{u})$ is negligible, we can try
directly the solution of the whole system by using techniques such as
the preconditioned CG (PCG) or the multigrid method
which in theory should suffer from the loss of symmetry in the linear
system.

In Figure \ref{fig:CG}a we plot the average (circles), minimum and
maximum (vertical lines) number of CG iterations performed during
the integration until final time, at different spatial
resolutions. A least square fit shows that the number of iterations
grows as $N^{0.5491}$. As previously observed in connection with the
GMRES method, the number of iterations is again essentially
proportional to $\sqrt{N}$, which agrees with the discussion in
Subsection \ref{ssec:MULTI}.

However, the number of iterations for fine grids is a lot higher
that the ones with GMRES (up to $950$ instead of $160$ with a grid
of $1024$ points) so the latter has to be preferred.

In a similar way, we consider $X_N(\vec{u})$ as preconditioner
in the PCG method. Results are shown in Figure \ref{fig:CG}b.  The
number of iteration is again essentially constant with respect to $N$,
but also in the preconditioned version, the GMRES is slightly better
since the number of PCG iterations, 8 or 9, is higher that the number
of PGMRES iterations which was equal to 6. Furthermore we have a
higher variance in the number of iterations, due to the weak
non-symmetry of whole matrix.

\subsubsection{MGM}
\label{sssec:num:MGM}
\begin{figure}
\hfil
\begin{tabular}{cc}
\includegraphics[width=.4\textwidth]{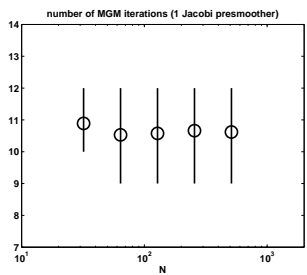}&
\end{tabular}
\hfil \caption{Average, minimum and maximum number of
  MGM iterations performed during the integration
  until final time.
  } \label{fig:MGM}
\end{figure}

We test the optimality of MGM, as discussed in Section
\ref{ssec:MULTI}.
We apply a single recursive call, that is the classical $V$-cycle
procedure. As smoother, we use a single Jacobi step with
damping factor equal to $2/3$. We observe mesh independent behaviour
with 10 or 11 iterations (see Figure \ref{fig:MGM}).

We also tried other more sophisticated multi-iterative approaches by
adding one step of post-smoother with Gauss-Seidel or standard Jacobi:
the number of iterations drops to 6, but the cost per iteration is
almost doubled, so that we do not observe a real advantage. The use of
one step of CG or one step of GMRES as post-smoother is not better
while one step of PGMRES with preconditioner equal to
$X_N(\vec{u})$ reduces the number of iterations, but not
enough compared with the cost of the solver needed for handling a
generic system with the preconditioner as coefficient matrix.

\subsubsection{Krylov methods with MGM as preconditioner}
\label{sssec:num:KrylovMGM}
\begin{figure}
  \centering{
    \begin{tabular}{cc}
      \includegraphics[width=0.4\textwidth]{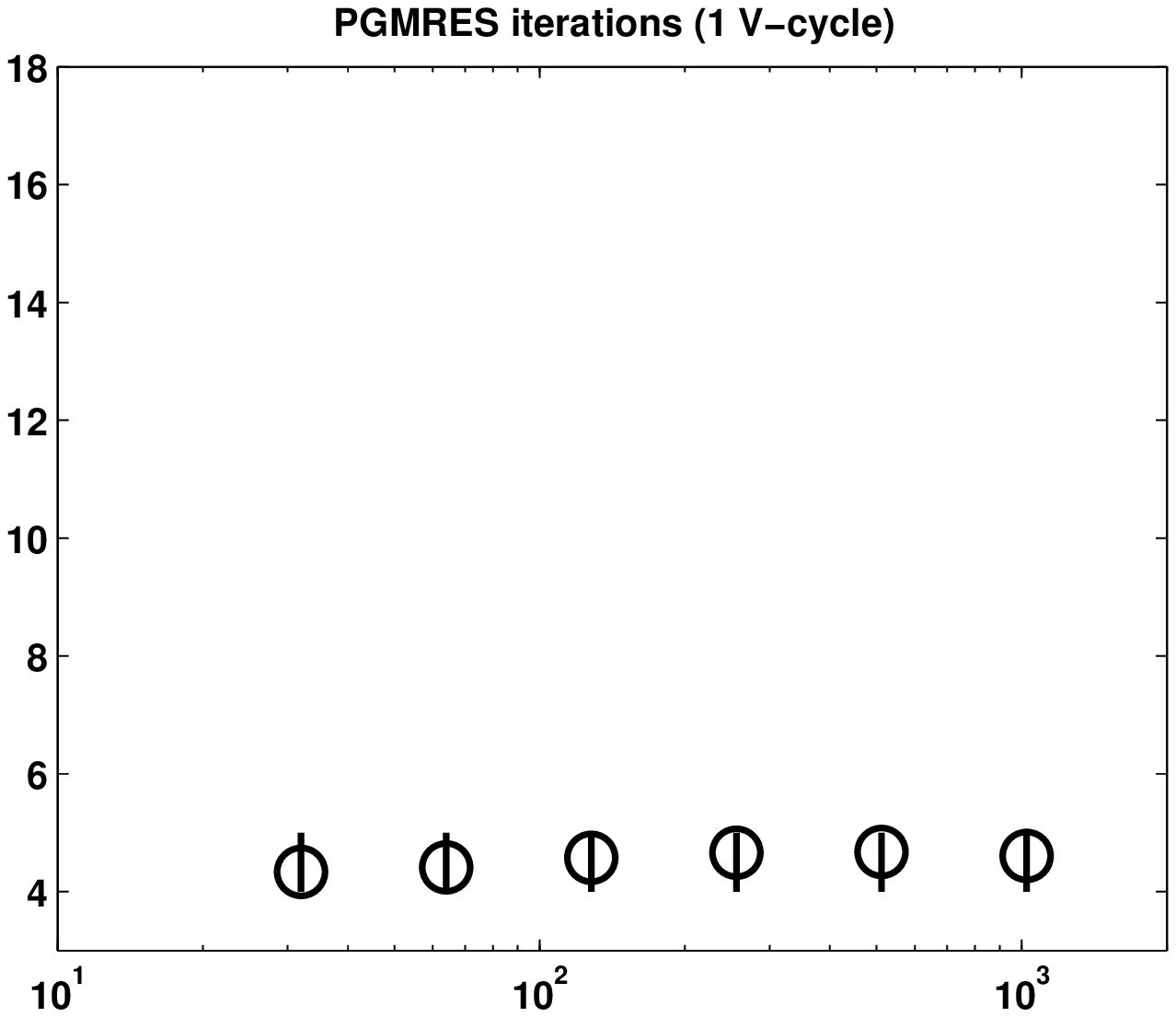}
      &
      \includegraphics[width=0.4\textwidth]{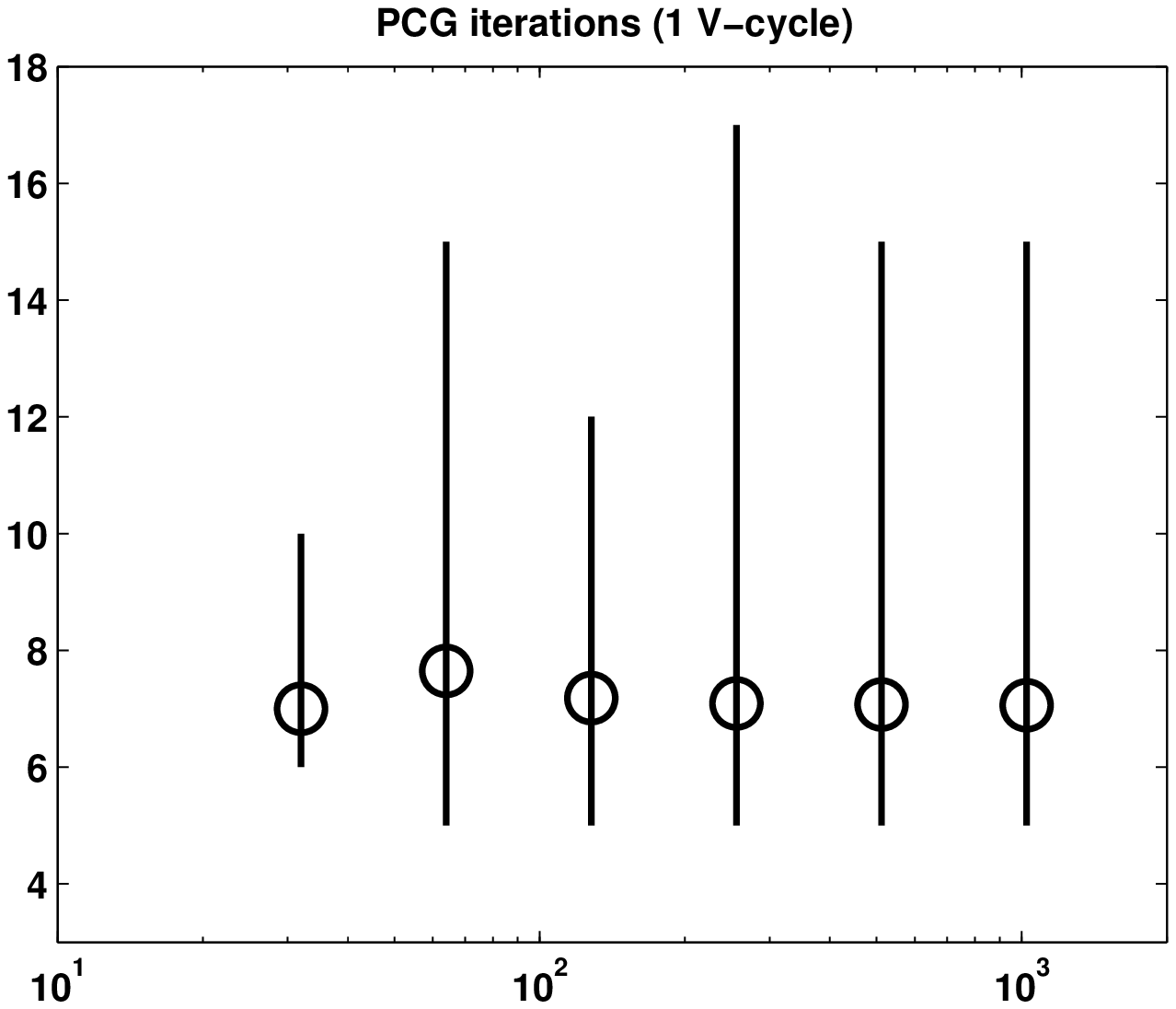}
      \\
      {\bf (a)} & {\bf (b)}
    \end{tabular}
  }
  \caption{Average, minimum and maximum number of PGMRES (a) and PCG
    (b) iterations performed during the integration until final time,
    when using one MGM V-cycle as preconditioner.
  }
  \label{fig:MGM+Krylov}
\end{figure}

The previous experiments confirm that the MGM is an excellent solver
for our linear system. Often this method is also applied as
preconditioner in a Krylov method instead of employing it as a
solver. In other words, as preconditioning step, we perform a
single V-cycle iteration, with the same coefficient matrix and
where the datum is the residual vector at the current iteration.

With the use of such very cheap MGM preconditioning, the PGMRES
converges within 4 or 5 iterations independently of the size of the
involved matrices (see Figure \ref{fig:MGM+Krylov}a). Comparing with
the GMRES method preconditioned with the symmetric part of $F^\prime(\vec{u})$
considered in \ref{sssec:num:GMRES} and Figure \ref{fig:GMRES}b, the
present preconditioning strategy is not only computationally cheaper,
but it is also more effective since it achieves a stronger reduction of
the number of GMRES iterations.

Analogously, the application of the same MGM preconditioning in the
PCG method leads to a convergence within 7 or 8 iterations, again
independently of the system sizes (see Figure
\ref{fig:MGM+Krylov}b).

In conclusion, V-cycle preconditioning in connection with
GMRES has to be preferred, taking into account the simplicity, the
robustness (less variance in the iteration count), and the number of
iterations. Indeed, due to the small iteration count, also the memory
requirement does not pose any difficulty, since the number of vectors
that have to be stored in the GMRES process is very reasonable.

\section{2D generalization}
\label{sec:2D}
In this section we describe a straightforward 2D generalization of the
numerical approach studied in the previous part of the paper. To this
end, we consider a rectangular domain
$\Omega=[a_0,a_1]\times[b_0,b_1]\subset\R^2$ and the grid points
$x_{i,j}=(a_0+ih,b_0+jk)$. For simplicity and without loss of
generality, we also assume that the region $\Omega$ is square and
choose identical discretization steps in the two directions
(i.e. $h=k$), so that using $N$ points per direction we have
$h=k=(a_1-a_0)/(N+1)=(b_1-b_0)/(N+1)$. The grid is thus composed of
the $(N+2)^2$ points $x_{i,j}$ for $i$ and $j$ ranging from $0$ to
$N+1$. We denote with $u_{i,j}$ the numerical value approximating
$u(x_{i,j})$. Of course, as in the one-dimensional case the use of
Dirichlet boundary conditions allows to reduce to gridding to the
$N^2$ internal points.

In this setting, we generalize the finite difference discretization
\eqref{eq:FDlaplacian} of the differential operator as follows:
\begin{multline}
\label{eq:diffop2d}
\left.
\frac{\partial}{\partial x} \left( D(u) \frac{\partial u}{\partial x}\right)
+\frac{\partial}{\partial y} \left( D(u) \frac{\partial u}{\partial y}\right)
\right|_{x=x_{i,j}}
\\=
\frac{D_{i+1/2,j}u_{i+1,j}%
  -\left(D_{i+1/2,j}+D_{i-1/2,j}\right)u_{i,j}%
  +D_{i-1/2,j}u_{i-1,j}%
}{h^2}
\\+
\frac{D_{i,j+1/2}u_{i,j+1}%
  -\left(D_{i,j+1/2}+D_{i,j-1/2}\right)u_{i,j}%
  +D_{i,j-1/2}u_{i,j-1}%
}{h^2}
+o(1),
\end{multline}
where we denoted
\[ D_{i+1/2,j} = \frac{D_{i+1,j}+D_{i,j}}{2}
\qquad
D_{i,j+1/2} = \frac{D_{i,j+1}+D_{i,j}}{2}.
\]
In order to write in matrix form the approximated differential
operator above, we must choose an ordering of the unknowns $u_{i,j}$,
arranging them into a vector $\vec{u}$ and approximate
\[ \left[\nabla\cdot(D(u)\nabla{u})(x_{i,j})\right]_{i,j=1}^n
\simeq \frac{1}{h^2}L_{D(\vec{u})}\,\vec{u}.
\]
The positions of the nonzero entries of the matrix $L_{D(\vec{u})}$ of
course depend on the chosen ordering, so here we keep a double-index
notation for the elements of $\vec{u}$ and of the matrix
entries. Therefore, following \eqref{eq:diffop2d}, $L_{D(\vec{u})}$
has entries
\begin{multline*}
\left[L_{D(\vec{u})}\right]_{i,j}^{l,m}
=
\delta_{i,l}\delta_{j,m}\left(
  -D_{i+1/2,j}-D_{i-1/2,j}-D_{i,j+1/2}-D_{i,j-1/2}
\right)
\\
+\delta_{l,i+1}\delta_{m,j}D_{i+1/2,j}
+\delta_{l,i-1}\delta_{m,j}D_{i-1/2,j}
\\
+\delta_{l,i}\delta_{m,j+1}D_{i,j+1/2}
+\delta_{l,i}\delta_{m,j-1}D_{i,j-1/2}
\end{multline*}
on the $(i,j)^{\text{th}}$ row and $(l,m)^{\text{th}}$ column. The
actual sparsity pattern of the resulting matrix thus depends on the
ordering of the unknowns $u_{i,j}$; with the usual lexicographic
ordering that has $u_{i,j}$ in the $(i+N(j-1))^{\text{th}}$ position of
the vector $\vec{u}$, one may have, as in the case of the standard
Laplacian operator, nonzero entries only on the main diagonal, on the
$1^{\text{st}}$ and $N^{\text{th}}$ upper and lower diagonals.

Each timestep with the Crandall Liggett formula \eqref{CLformula} thus
requires the solution of the nonlinear equation defined by
\[ F(\vec{u})
= \vec{u}
- \frac{\DT}{h^2} L_{D(\vec{u})}\,\vec{u}
- \vec{u}^{n-1}.
\]
As in the one-dimensional case we propose to approximate the solution
of the nonlinear equation with the Newton's method; an analysis similar to that of
Theorem \ref{prop:newton} can be carried out in the new 2D context.
The Jacobian of $F(\vec{u})$ is
\begin{equation} \label{Jf2d}
F^\prime(\vec{u}) = I - \frac{\DT}{h^2}L_{D(\vec{u})} -\frac{\DT}{h^2}Y(\vec{u})
\end{equation}
where
\begin{equation}
Y_{i,j}^{l,m}(\vec{u}) = \sum_{l,m}
\frac{\partial
  \left[L_{D((\vec{u}))}\right]_{i,j}^{l,m}}{\partial u_{r,s}}
u_{l,m}
\end{equation}
with the double-index notation as above.

A tedious but straightforward computation yields
\begin{multline*}
Y_{i,j}^{r,s}(\vec{u}) =
\frac12 D^\prime_{i,j}\delta_{i,r}\delta_{j,s}
\left(-4u_{i,j}+u_{i+1,j}+u_{i-1,j}+u_{i,j+1}+u_{i,j-1}\right)
\\
+\frac12D^\prime_{i+1,j}\delta_{r,i+1}\delta_{s,j}\left(u_{i+1,j}-u_{i,j}\right)
+\frac12D^\prime_{i-1,j}\delta_{r,i-1}\delta_{s,j}\left(u_{i-1,j}-u_{i,j}\right)
\\
+\frac12D^\prime_{i,j+1}\delta_{r,i}\delta_{s,j+1}\left(u_{i,j+1}-u_{i,j}\right)
+\frac12D^\prime_{i,j-1}\delta_{r,i}\delta_{s,j-1}\left(u_{i,j-1}-u_{i,j}\right)
\end{multline*}
for the generic entry of $Y(\vec{u})$. (The obvious changes must be
taken into account to implement the boundary conditions, e.g. either
eliminating the unknowns for the points on the Dirichlet boundary or
the unknowns on suitably chosen ghost points outside the Neumann
boundary.)

\begin{figure}
  \begin{tabular}{cc}
    \includegraphics[width=0.45\linewidth]{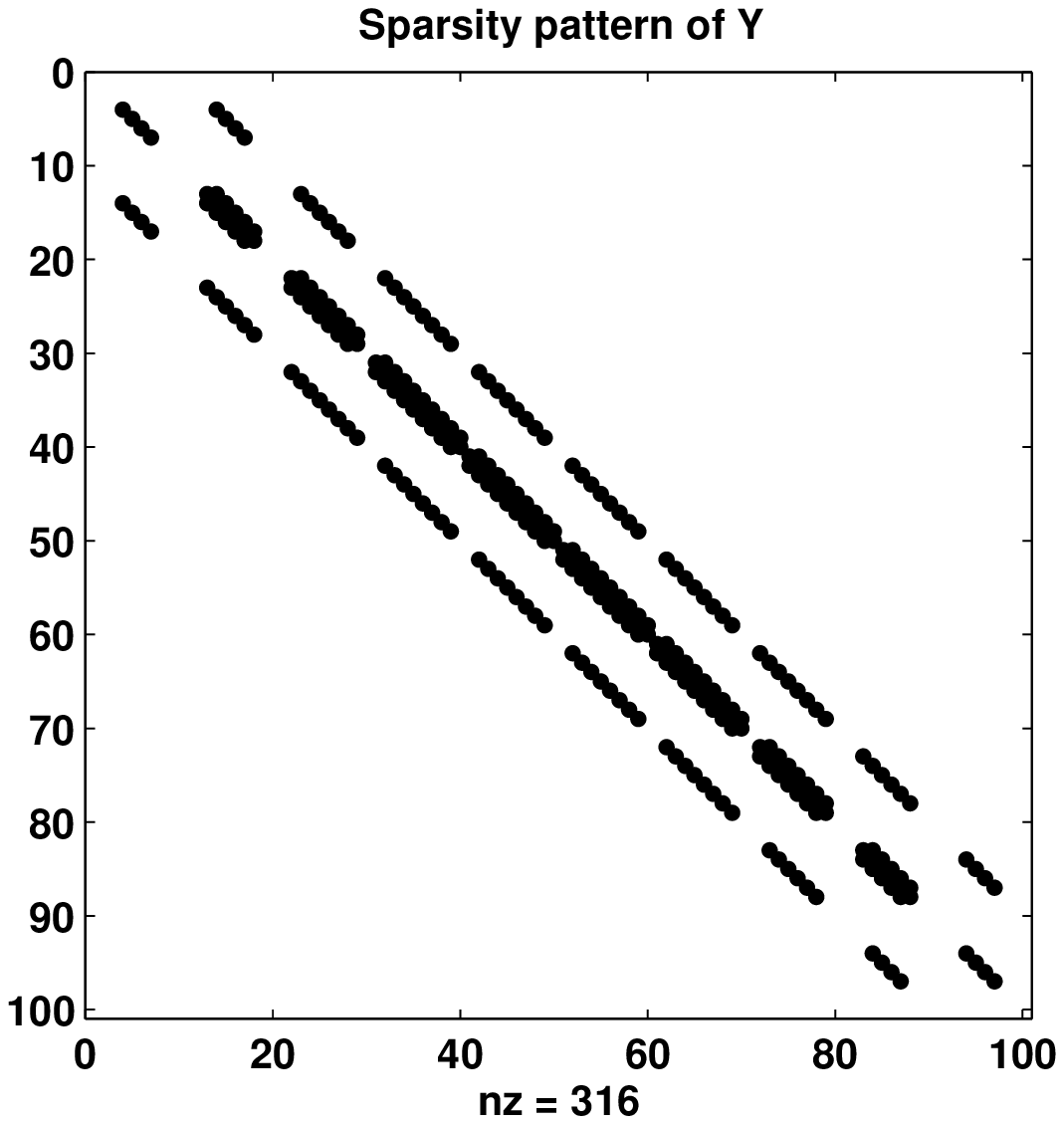}
    &
    \includegraphics[width=0.45\linewidth]{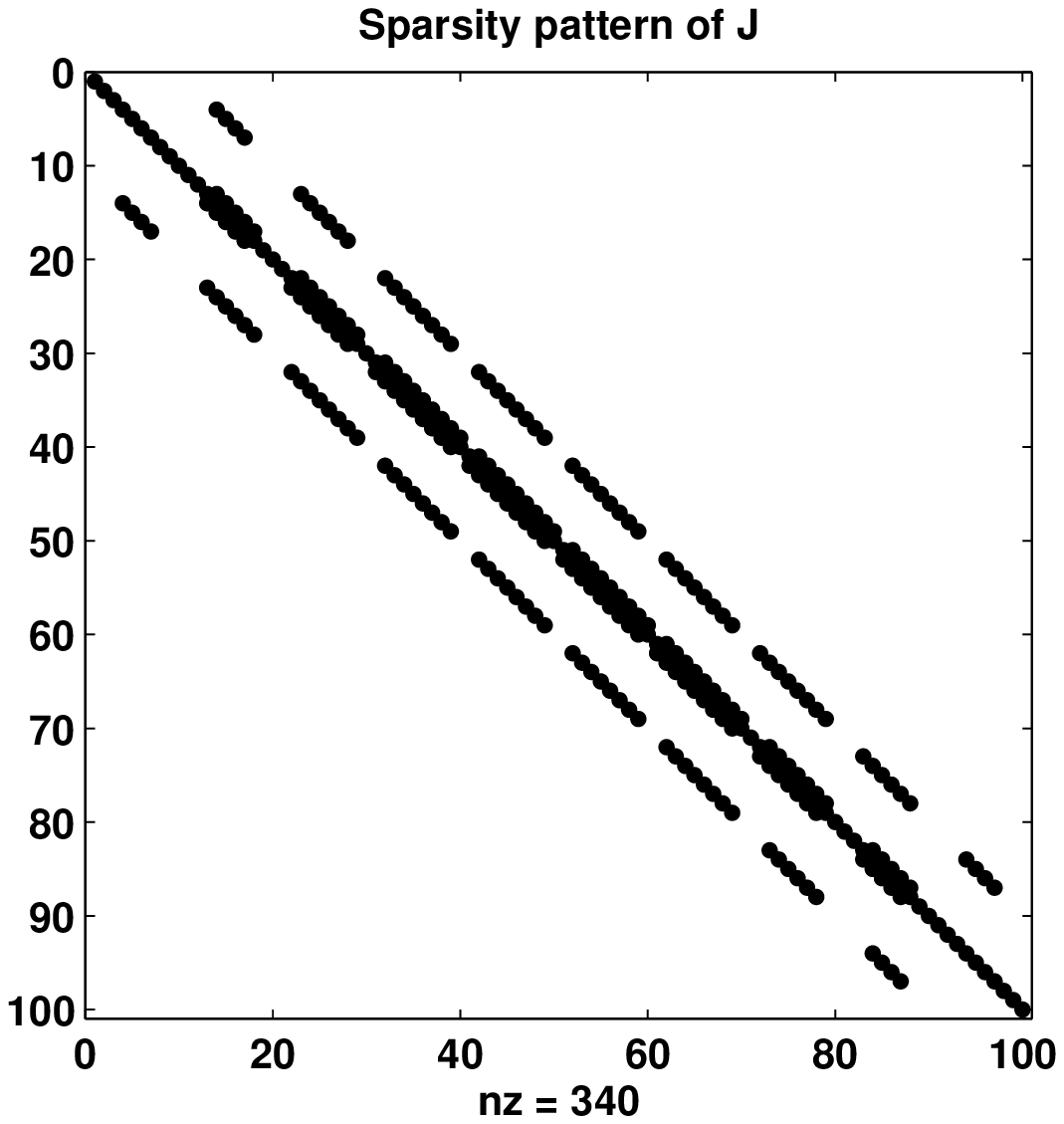}\\
    {\bf (a)}&{\bf (b)}
  \end{tabular}
  \caption{Sparsity pattern of $Y(\vec{u})$ (a) and
    $F^\prime(\vec{u})$ (b) on a $10\times10$ grid with the unknowns in
    lexicographic ordering.}
  \label{fig:sparsy2d}
\end{figure}

As in the one-dimensional case, the matrix $Y(\vec{u})$ may be written
as $Y(\vec{u})=T(\vec{u})D^\prime(\vec{u})$ where $D^\prime(\vec{u})$
is a diagonal matrix with entries equal to $D^\prime(u_{i,j})$ and in
smooth regions of the solution, the nonzero entries of $Y(\vec{u})$
are $O(h^2)$ on the main diagonal and $O(h)$ outside. Moreover,
entries of $Y(\vec{u})$ expected to be nonzero may in fact be null
because the approximate solution is locally flat in a neighbourhood or
because some of the $D^{\prime}(u_k)$ may be null.  When using the
natural ordering of the unknowns described above, the sparsity
patterns of $Y(\vec{u})$ and $F^\prime(\vec{u})$ for the Barenblatt
solution are illustrated in figure \ref{fig:sparsy2d}. The gaps along
the diagonals of $Y(\vec{u})$ correspond to the regions where the
approximate solution $\vec{u}$ is flat.

\begin{figure}
  \centering\includegraphics[width=0.4\linewidth]{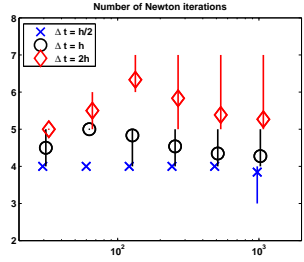}
  \caption{Average, minimum and maximum number of Newton iterations
    performed during the integration until final time. The $3$ data
    series for each $N$ have been slightly shifted for clarity.}
  \label{fig:Newton2d}
\end{figure}

We performed our tests with the two-dimensional Barenblatt solution
\citep[see][]{VazquezBOOK} with exponent $m=4$ on grids of size
$N\times{N}$ for $N$ ranging from $32$ to $1024$.
First of all we note that the number of Newton iterations required at
each timestep is almost independent of $N$ and is (on average) $4$
when $\DT=0.5h$, $4.5$ when $\DT=h$ and $6.5$ when $\DT=2h$ (see
Figure \ref{fig:Newton2d}).

\begin{figure}
  \begin{tabular}{cc}
    \includegraphics[width=0.4\linewidth]{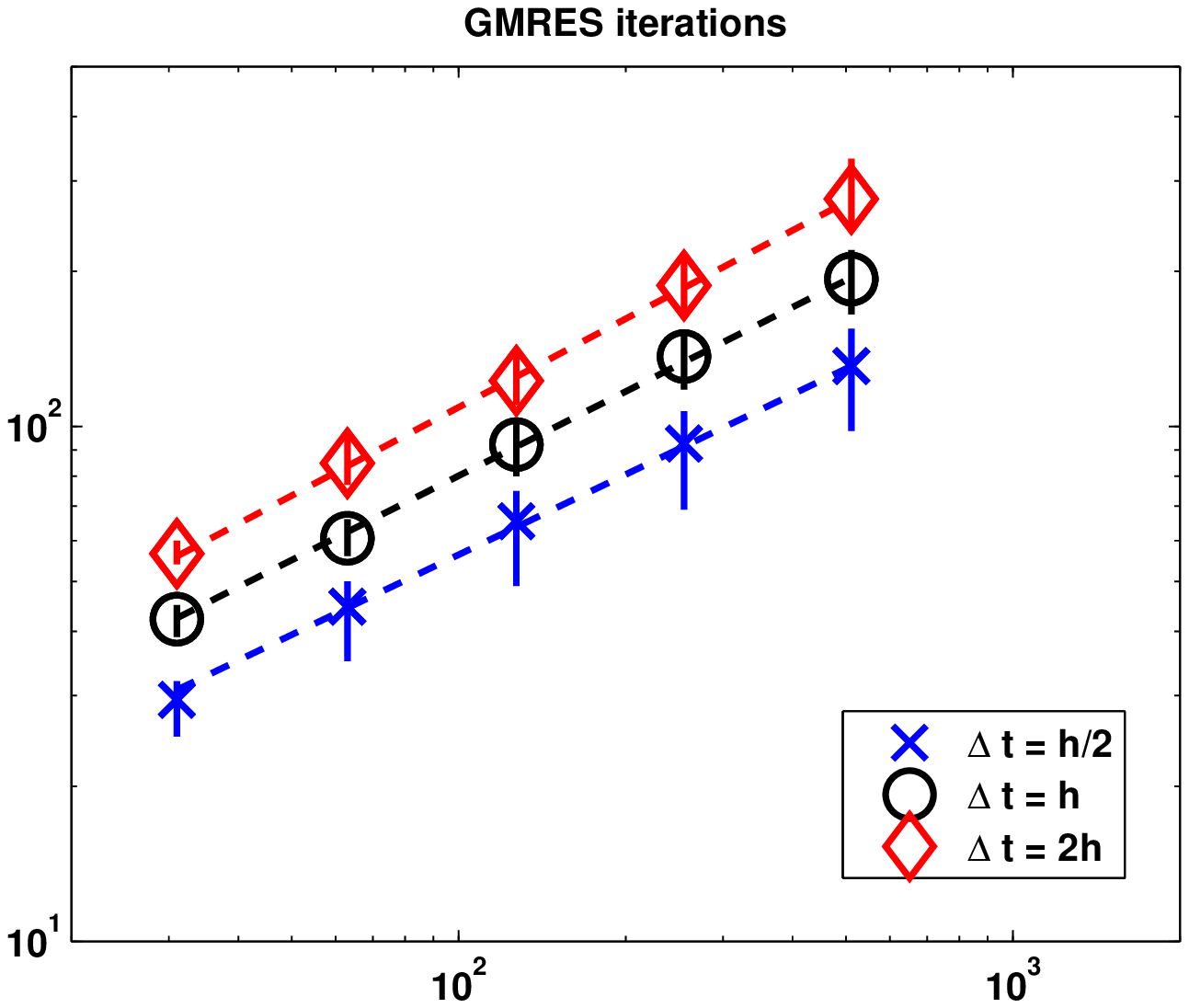}
    &
    \includegraphics[width=0.4\linewidth]{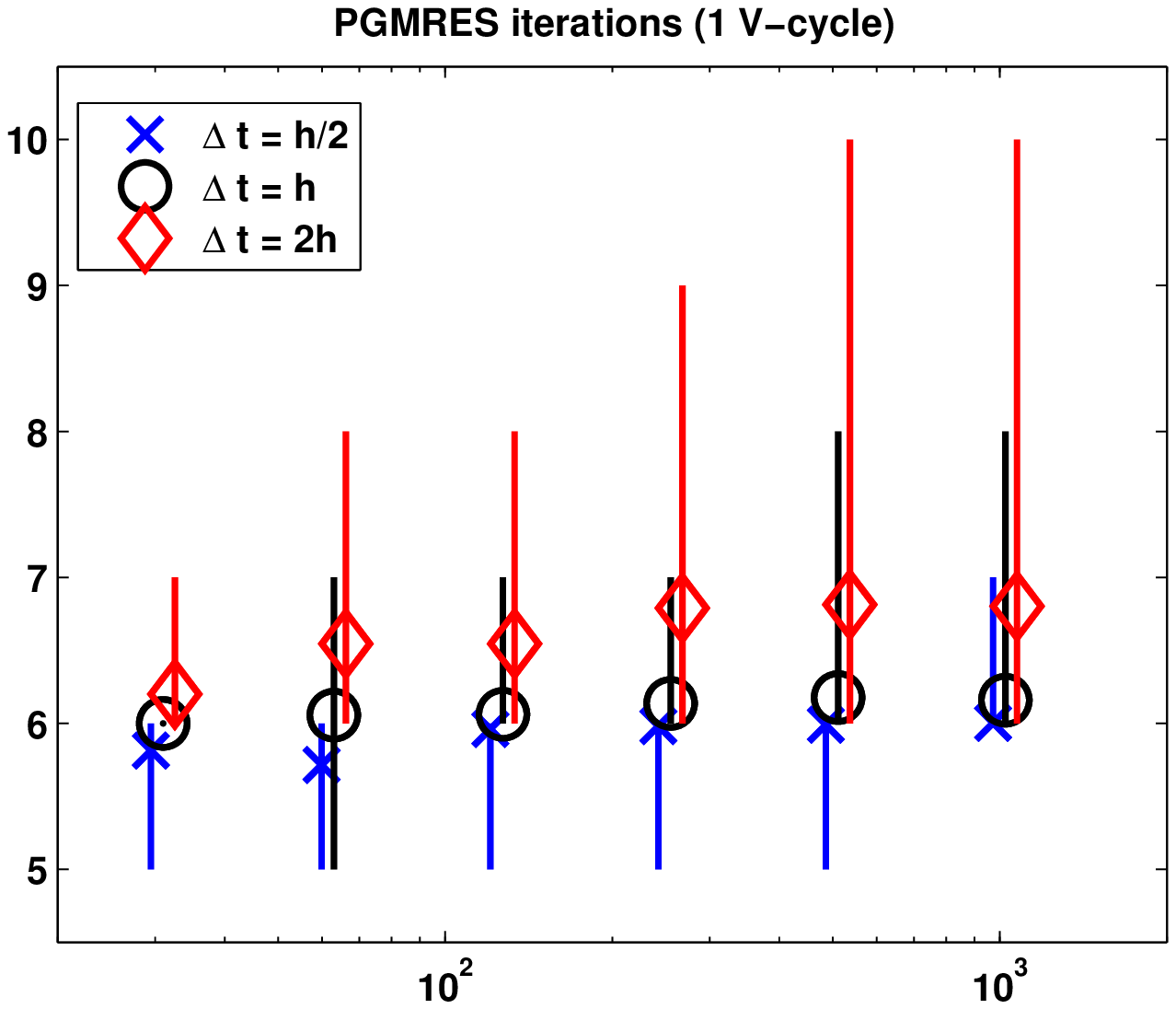}
    \\
    (a)&(b)
  \end{tabular}
  \caption{Number of GMRES iterations at different grid sizes, in
    2D. (a) without preconditioning. (b) with V-cycle
    preconditioner. On the right, the $3$ data series for each $N$
    have been slightly shifted for clarity.}
  \label{fig:GMRES2d}
\end{figure}

We point out that the results of Section \ref{ssec:SPECTRAL}
generalize to the two-dimensional case and thus we perform numerical
tests using a multigrid iteration as preconditioner for PGMRES
which in Section \ref{sec:NUM} provided best results in the one dimensional case.

In particular we employ a single V-cycle iteration, with a
Galerkin approach using the bilinear interpolation as prolongation
operator and one step of red-black Gauss-Seidel as pre-smoother.  In
Figure \ref{fig:GMRES2d} we plot the mean (symbols) and
minimum-maximum (solid lines) number of GMRES iterations needed at
different spatial resolutions.
Different colours correspond to different choices of $\DT$, namely $\DT=h/2$ (blue
crosses), $\DT=h$ (black circles) and $\DT=2h$ (red diamonds). The
left panel shows that, without preconditioning, the number of GMRES
iterations grows with the grid size: least square fits yield the
approximations $N^{0.5165}$, $N^{0.5435}$ and $N^{0.5702}$
respectively for the number of GMRES iterations on an $N\times{N}$
grid with the three choices of $\DT$ mentioned above. For homogeneity,
the results for $N=1024$ are not reported in the graph, since they
require the restarted GMRES method or a parallel implementation, due
to memory limitations when run on a PC with 8Mb of RAM.

Figure \ref{fig:GMRES2d}b clearly demonstrates the optimality of the
preconditioning strategy adopted, with the number of iterations
being in the narrow range $5$--$10$ when $N$ ranges from $32$ to
$1024$ and with all the three choices of the time step and with the
average number of iterations being always between $5$ and $7$. We
note in passing that we also employed damped Jacobi as a smoother
with analogous results on the optimality, but observing a slightly
higher number of iterations ($8$--$11$ on average).


\section{Conclusions and future developments}
\label{sec:CONCL}

The novel contribution of this paper relies in the proposal of a fully
implicit numerical method for dealing with nonlinear degenerate
parabolic equations, in its convergence and stability analysis, and in
the study of the related computational cost. Indeed the nonlinear
nature of the underlying mathematical model requires the application
of a fixed point scheme.  We identified the classical Newton
method in which, at every step, the solution of a large, locally
structured, linear system has been handled by using specialised
iterative or multi-iterative solvers. In particular, we provide
a spectral analysis of the relevant matrices which has been crucial
for identifying appropriate preconditioned Krylov methods with
efficient V-cycle preconditioners.  Numerical experiments for the
validation of our multi-facet analysis complement this
contribution.

Among the vast range of possible applications of degenerate parabolic
equations, we point out a recent one in the field of monument
conservation in \citet{SDS:monum}, where an approximation technique
derived from the one analysed here has been successfully employed in
the forecast of marble deterioration on monuments.
Having in mind the application to more complicated monument geometry,
we will pursue the extension of the results of this paper to the case
of finite element methods for the space discretization.

\bibliographystyle{IMANUM-BIB}
\bibliography{DegDiffMGM}

\begin{thebibliography}{}

\bibitem[Al{\`i} {\em et~al.}(2007)Al{\`i}, Furuholt, Natalini, \&
  Torcicollo]{AFNT07:model}
{\sc Al{\`i}, G., Furuholt, V., Natalini, R. \& Torcicollo, I.} (2007)
\newblock A mathematical model of sulphite chemical aggression of limestones
  with high permeability. {I}. {M}odeling and qualitative analysis.
\newblock {\em Transp. Porous Media\/}, {\bf 69}, 109--122.

\bibitem[Aregba~Driollet {\em et~al.}(2004)Aregba~Driollet, Diele, \&
  Natalini]{ADN:sulfation}
{\sc Aregba~Driollet, D., Diele, F. \& Natalini, R.} (2004)
\newblock A mathematical model for the {$\mathrm{SO}_2$} aggression to calcium
  carbonate stones: numerical approximation and asymptotic analysis.
\newblock {\em SIAM J. Appl. Math.}, {\bf 64}, 1636--1667.

\bibitem[Beckermann \& Serra-Capizzano(2007)Beckermann \&
  Serra-Capizzano]{BS07}
{\sc Beckermann, B. \& Serra-Capizzano, S.} (2007)
\newblock On the asymptotic spectrum of finite element matrix sequences.
\newblock {\em SIAM J. Numer. Anal.}, {\bf 45}, 746--769 (electronic).

\bibitem[Berger {\em et~al.}(1979)Berger, Brezis, \& Rogers]{BBR79}
{\sc Berger, A., Brezis, H. \& Rogers, J.} (1979)
\newblock A numerical method for solving the problem {$u_t-\Delta f(u)=0$}.
\newblock {\em RAIRO numerical analysis\/}, {\bf 13}, 297--312.

\bibitem[Bertaccini {\em et~al.}(2005)Bertaccini, Golub, Serra~Capizzano, \&
  Tablino~Possio]{BGST05}
{\sc Bertaccini, D., Golub, G.~H., Serra~Capizzano, S. \& Tablino~Possio, C.}
  (2005)
\newblock Preconditioned {HSS} methods for the solution of non-{H}ermitian
  positive definite linear systems and applications to the discrete
  convection-diffusion equation.
\newblock {\em Numer. Math.}, {\bf 99}, 441--484.

\bibitem[Bhatia(1997)Bhatia]{Bhatia:book}
{\sc Bhatia, R.} (1997)
\newblock {\em Matrix analysis\/}. Graduate Texts in Mathematics,  vol. 169.
\newblock New York: Springer-Verlag, pp. xii+347.

\bibitem[Br{\'e}zis \& Pazy(1972)Br{\'e}zis \& Pazy]{BP72}
{\sc Br{\'e}zis, H. \& Pazy, A.} (1972)
\newblock Convergence and approximation of semigroups of nonlinear operators in
  {B}anach spaces.
\newblock {\em J. Functional Analysis\/}, {\bf 9}, 63--74.

\bibitem[Cavalli {\em et~al.}(2007)Cavalli, Naldi, Puppo, \&
  Semplice]{CNPS07:degdiff}
{\sc Cavalli, F., Naldi, G., Puppo, G. \& Semplice, M.} (2007)
\newblock High-order relaxation schemes for non linear degenerate diffusion
  problems.
\newblock {\em SIAM Journal on Numerical Analysis\/}, {\bf 45}, 2098--2119.

\bibitem[Clarelli {\em et~al.}(2009)Clarelli, Giavarini, Natalini, Nitsch, \&
  Santarelli]{CGNNS:teos}
{\sc Clarelli, F., Giavarini, C., Natalini, R., Nitsch, C. \& Santarelli, M.}
  (2009)
\newblock Mathematical models for the consolidation processes in stones.
\newblock {\em Proc. of ``{I}nternational {S}ymposium: {S}tone {C}onsolidation
  in {C}ultural {H}eritage - research and practice''. {L}isbona, {M}ay 2008.}
\newblock
\newblock to appear.

\bibitem[Crandall \& Liggett(1971)Crandall \& Liggett]{CL71}
{\sc Crandall, M. \& Liggett, T.} (1971)
\newblock Generation of {S}emi-{G}roups of non linear transformations on
  general {B}anach spaces.
\newblock {\em Amer. J. Math.}, {\bf 93}, 265--298.

\bibitem[Golinskii \& Serra-Capizzano(2007)Golinskii \&
  Serra-Capizzano]{GS07:jacobi-nonsymm}
{\sc Golinskii, L. \& Serra-Capizzano, S.} (2007)
\newblock The asymptotic properties of the spectrum of nonsymmetrically
  perturbed {J}acobi matrix sequences.
\newblock {\em J. Approx. Theory\/}, {\bf 144}, 84--102.

\bibitem[Golub \& Van~Loan(1996)Golub \& Van~Loan]{Golub:book}
{\sc Golub, G.~H. \& Van~Loan, C.~F.} (1996)
\newblock {\em Matrix computations\/}.
\newblock Johns Hopkins Studies in the Mathematical Sciences, third edn.
\newblock Baltimore, MD: Johns Hopkins University Press, pp. xxx+698.

\bibitem[Greenbaum(1997)Greenbaum]{Greenbaum:book}
{\sc Greenbaum, A.} (1997)
\newblock {\em Iterative methods for solving linear systems\/}. Frontiers in
  Applied Mathematics,  vol.~17.
\newblock Philadelphia, PA: Society for Industrial and Applied Mathematics
  (SIAM), pp. xiv+220.

\bibitem[Hackbusch(1985)Hackbusch]{Hackbusch:book}
{\sc Hackbusch, W.} (1985)
\newblock {\em Multigrid methods and applications\/}. Springer Series in
  Computational Mathematics,  vol.~4.
\newblock Berlin: Springer-Verlag, pp. xiv+377.

\bibitem[Magenes {\em et~al.}(1987)Magenes, Nochetto, \& Verdi]{MNV87}
{\sc Magenes, E., Nochetto, R.~H. \& Verdi, C.} (1987)
\newblock Energy error estimates for a linear scheme to approximate nonlinear
  parabolic problems.
\newblock {\em RAIRO Mod{\'e}l. Math. Anal. Num{\'e}r.}, {\bf 21}, 655--678.

\bibitem[Ortega \& Rheinboldt(1970)Ortega \& Rheinboldt]{Ortega:book}
{\sc Ortega, J.~M. \& Rheinboldt, W.~C.} (1970)
\newblock {\em Iterative solution of nonlinear equations in several
  variables\/}.
\newblock New York: Academic Press, pp. xx+572.

\bibitem[Saad(2003)Saad]{Saad:book}
{\sc Saad, Y.} (2003)
\newblock {\em Iterative methods for sparse linear systems\/}, second edn.
\newblock Philadelphia, PA: Society for Industrial and Applied Mathematics, pp.
  xviii+528.

\bibitem[Semplice {\em et~al.}(2009)Semplice, Donatelli, \&
  Serra-Capizzano]{SDS:monum}
{\sc Semplice, M., Donatelli, M. \& Serra-Capizzano, S.} (2009)
\newblock Preconditioned fully implicit pde solvers for degenerate parabolic
  equations with applications to monument conservation.
\newblock {\em http:{$\backslash\backslash$}www.arXiv.org\/}, {\bf
  0907.2600v1}.

\bibitem[Serra-Capizzano(1993)Serra-Capizzano]{Serra93:multi}
{\sc Serra-Capizzano, S.} (1993)
\newblock Multi-iterative methods.
\newblock {\em Comput. Math. Appl.}, {\bf 26}, 65--87.

\bibitem[Serra-Capizzano(2006)Serra-Capizzano]{Serra06:glt}
{\sc Serra-Capizzano, S.} (2006)
\newblock The {GLT} class as a generalized {F}ourier analysis and applications.
\newblock {\em Linear Algebra Appl.}, {\bf 419}, 180--233.

\bibitem[Serra-Capizzano \& Tablino-Possio(2004)Serra-Capizzano \&
  Tablino-Possio]{ST04}
{\sc Serra-Capizzano, S. \& Tablino-Possio, C.} (2004)
\newblock Multigrid methods for multilevel circulant matrices.
\newblock {\em SIAM J. Sci. Comput.}, {\bf 26}, 55--85 (electronic).

\bibitem[Stoer \& Bulirsch(2002)Stoer \& Bulirsch]{Sto02}
{\sc Stoer, J. \& Bulirsch, R.} (2002)
\newblock {\em Introduction to numerical analysis\/}. Texts in Applied
  Mathematics,  vol.~12, third edn.
\newblock New York: Springer-Verlag, pp. xvi+744.
\newblock Translated from the German by R. Bartels, W. Gautschi and C.
  Witzgall.

\bibitem[Tilli(1998)Tilli]{Tilli98}
{\sc Tilli, P.} (1998)
\newblock Locally {T}oeplitz sequences: spectral properties and applications.
\newblock {\em Linear Algebra Appl.}, {\bf 278}, 91--120.

\bibitem[Trottenberg {\em et~al.}(2001)Trottenberg, Oosterlee, \&
  Sch{\"u}ller]{Trottenberg:book}
{\sc Trottenberg, U., Oosterlee, C.~W. \& Sch{\"u}ller, A.} (2001)
\newblock {\em Multigrid\/}.
\newblock San Diego, CA: Academic Press Inc., pp. xvi+631.
\newblock With contributions by A. Brandt, P. Oswald and K. St{\"u}ben.

\bibitem[Varga(1962)Varga]{Varga:matrixbook}
{\sc Varga, R.~S.} (1962)
\newblock {\em Matrix iterative analysis\/}.
\newblock Englewood Cliffs, N.J.: Prentice-Hall Inc., pp. xiii+322.

\bibitem[V\'azquez(2007)V\'azquez]{VazquezBOOK}
{\sc V\'azquez, J.~L.} (2007)
\newblock {\em The porous medium equation\/}.
\newblock Oxford Mathematical Monographs.
\newblock Oxford: The Clarendon Press Oxford University Press, pp. xxii+624.
\newblock Mathematical theory.

\end{thebibliography}

\end{document}